\documentclass[10pt]{amsart}

\usepackage{amsthm,hyperref,amsmath,amssymb,mathtools}
\usepackage{graphicx,subfig,color,indentfirst,epstopdf,booktabs}

\usepackage[font=small]{caption}

\numberwithin{equation}{section}

\theoremstyle{definition}
\newtheorem{definition}{Definition}

\hypersetup{
	colorlinks=true,
	allcolors=blue
}

\begin{document}
\baselineskip=14pt

\title[Numerical analysis of the Gatenby-Gawlinski model reductions]{Numerical investigation\\ of some reductions for\\ the Gatenby-Gawlinski model}

\author[C. Mascia]{Corrado Mascia}
\address[Corrado Mascia]{Dipartimento di Matematica G. Castelnuovo\\ Sapienza Universit\`a di Roma\\ Piazzale Aldo Moro 2 - 00185 Roma (Italy)}
\email{corrado.mascia@uniroma1.it}

\author[P. Moschetta]{Pierfrancesco Moschetta}
\address[Pierfrancesco Moschetta]{Dipartimento di Matematica G. Castelnuovo\\ Sapienza Universit\`a di Roma\\ Piazzale Aldo Moro 2 - 00185 Roma (Italy)}
\email{moschetta@mat.uniroma1.it}

\author[C. Simeoni]{Chiara Simeoni}
\address[Chiara Simeoni]{Laboratoire de Math\'ematiques J.A. Dieudonn\'e CNRS UMR 7351, Universit\'e C\^ote D'Azur, Parc Valrose - 06108 Nice Cedex 2 (France)}
\email{chiara.simeoni@univ-cotedazur.fr}

\begin{abstract}
The Gatenby-Gawlinski model for cancer invasion is object of analysis in order to investigate the mathematical framework behind the model working by means of suitable reductions. We perform numerical simulations to study the sharpness/smoothness of the traveling fronts starting from a brief overview about the full model and proceed by examining the case of a two-equations-based and one-equation-based reduction. We exploit a numerical strategy depending on a finite volume approximation and employ a space-averaged wave speed estimate to quantitatively approach the traveling waves phenomenon. Concerning the one-equation-based model, we propose a reduction framed within the degenerate reaction-diffusion equations field, which proves to be effective in order to qualitatively recover the typical trends arising from the Gatenby-Gawlinski model. Finally, we carry out some numerical tests in a specific case where the analytical solution is available.
\end{abstract}



\maketitle

%%%------------------------------------------------------------------------------------------

\section{Introduction}
\label{sec:introduction}
The contribution of mathematical modeling as an effective investigation tool in the biomedical field, is becoming more and more relevant to the present time. Cancer research provides doubtless several interesting research paths, whose mathematical exploration is currently taking place on the back of promising results~\cite{Gatenby 2002}. In this scenario, the so-called \emph{Warburg effect}~\cite{Warburg 1930, Warburg 1956} plays a significant role within the phenomenological framework which is worth being carefully examined in order to accomplish a better understanding of the dynamics tumour growth is ruled by. In the 1920s, Otto Warburg~\cite{Warburg 1930} experimentally noticed that cancerous cells essentially rely on glycolytic metabolism, regardless of the oxygen availability: as a matter of fact, oxygen turns out to be the principle resource to allow normal cells performing glucose metabolism, due to the best yield of adenosine triphosphate (ATP) production by using oxidative phosphorylation. On the other hand, tumour cells appear to lean more towards glycolysis, leading to lactic acid fermentation. This is crucial for the so-termed acid-mediated invasion hypothesis~\cite{Bertuzzi 2010, Bertuzzi 2009, Gatenby 2003, Gatenby 2004, Gatenby 2006, Gatenby 2007, Smallbone 2008}, whose key point consists in assuming that acidification induced by lactic acid sets up a toxic microenvironment for normal cells and favors cancer cells spreading.

From a mathematical point view, all these qualitative statements are properly framed by the Gatenby-Gawlinski reaction-diffusion model~\cite{Gatenby 1996}, whose investigation has been being carried out by both numerical and analytical approaches~\cite{Davis 2018, Fasano 2009, Gatenby 1996, McGillen 2014, Moschetta 2019}. Basically, the model is developed for describing tumour cells proliferation at the expense of the local healthy tissue, exploiting the framework defined by the species evolution and assuming that carcinogenesis has already been performed. Indeed, the focal point is the interaction involving cancerous and healthy species at the tumour-host interface, through the lactic acid mediation.

The main purpose of this article lies in allowing a better understanding of the mathematical features the model is characterized by, specifically employing suitable hypothesis for building some model reductions, on the heels of what has already introduced in~\cite{Moschetta 2019} for the one-dimensional case. By means of numerical simulations, the traveling waves phenomenon is analyzed with emphasis on the qualitative structure of the fronts.
%which clearly exhibit a typical sharp trend~\cite{Sanchez 1995} for the reductions case.    

The contents of this paper are organized as follows. In Section~\ref{sec:twoeq}, a brief report about the general form of the model and information about its working are provided, along with the previous achievements described in~\cite{Moschetta 2019} concerning the two-equations reduction; afterwards, an investigation on the fronts is carried out by means of numerical simulations and a sensitivity analysis with respect to some system parameters is performed as well, considering as unknown the wave speed of the front. For this goal, a suitable space-averaged wave speed approximation is taken into account~\cite{LeVeque 1990}. Section~\ref{sec:oneeq} is aimed at building a one-equation-based reduction theoretically framed in the degenerate reaction-diffusion equations field~\cite{Malaguti 2003, Sanchez 1994 first, Sanchez 1994 second, Sanchez 1995, Sanchez 1996, Sanchez 1997}; the degenerate diffusion arising from the model reduction is almost everywhere differentiable in $[0,1]$ so that the sharpness of the fronts is to be checked by means of numerical simulations, due to the requirement for more regularity needed by the corresponding analytical results~\cite{Malaguti 2003, Sanchez 1994 first, Sanchez 1995, Sanchez 1996}. Moreover, numerical simulations are performed in order to qualitatively retrieve the dynamics exhibited by the Gatenby-Gawlinski model and several tests are carried out involving the related exact solution when its availability is ensured. Finally, in Section~\ref{sec:theend}, we discuss the conclusions of the manuscript and perspectives on future research.

%%%------------------------------------------------------------------------------------------

\section{Two-equations-based model reduction}
\label{sec:twoeq}

\subsection{Derivation and previous achievements}
The original Gatenby-Gawlinski model~\cite{Gatenby 1996} is composed by three equations (two PDEs plus one ODE) that, for the sake of convenience, we choose to make non-dimensionalized~\cite{Gatenby 1996, Moschetta 2019} (see~\cite{McGillen 2014} for a generalized version), so that we get 
\begin{equation}
\label{eqn:system}
\begin{dcases}
\frac{\partial u}{\partial t} = u (1-u) - d u w\\
\frac{\partial v}{\partial t} = r v (1-v ) + D \frac{\partial}{\partial x} \left[(1-u) \frac{\partial v}{\partial x}\right]\\
\frac{\partial w}{\partial t} = c (v-w) + \frac{\partial^{2} w}{\partial x^{2}}
\end{dcases}
\end{equation}
where the interval $[-1,1]$ is the one-dimensional domain, with $t \geq 0$, while $u(x,t)$, $v(x,t)$ and $w(x,t)$ are the unknown scaled functions which stand for the healthy tissue density, the tumour tissue density and the extracellular lactic acid concentration in excess, respectively. The densities $u$ and $v$ follow a logistic growth with normalized carrying capacities; $d$ is a death rate proportional to $w$, considered for reproducing healthy cells degradation brought on by lactic acid, while $r$ is a growth rate. Talking about the second equation, it is remarkable to recognize the structure of the degenerate diffusion term, in which $D$ is the diffusion constant for cancerous cells when the healthy tissue has already been degraded; on the contrary, when the local healthy cells concentration is equal to its normalized ceiling, the tumour cannot spread out as a consequence of a defense process of confinement~\cite{Gatenby 1996, Sanchez 1995}. Finally, in the third equation, the parameter $c$ plays the role of both a growth rate for the acid production (proportional to $v$) and a physiological reabsorption rate. For the boundary conditions, the homogeneous Neumann problem is set out.

As concerns the dynamics provided by~\eqref{eqn:system}, the results can be essentially condensed through two different kinds of behaviours~\cite{Gatenby 1996, McGillen 2014, Moschetta 2019}, both being framed within the propagating fronts theory: the first one, which happens in the regime $d<1$, is called \textit{heterogeneous invasion}, because of the coexistence of tumour and healthy tissues behind the wave front; if $d>1$, instead, we face a more aggressive invasion, the so-called \textit{homogeneous invasion}, due to the complete destruction of the healthy tissue out by the cancerous cells wave front. In this last regime, we recall the presence of a \textit{tumour-host hypocellular interstitial gap\,}, namely a separation zone between the healthy and cancer cells densities.

Information about the wave speed is object of analysis as well, showing that, in agreement with its more aggressive nature, the homogeneous invasion turns out to happen faster than the heterogeneous one~\cite{Gatenby 1996, McGillen 2014, Moschetta 2019}.

As far as the possibility of relying on a simplified version of the model \eqref{eqn:system}, we take advantage of what has been proposed in~\cite{Moschetta 2019} and afterwards try to go further. The assumption allowing to get a two-equations-based reduction is $w=v$, which leads to 
\begin{equation}
\label{eqn:system2}
\begin{dcases}
\frac{\partial u}{\partial t} = u (1-u) - d u v\\
\frac{\partial v}{\partial t} = r v (1-v) + D \frac{\partial}{\partial x} \left[(1-u) \frac{\partial v}{\partial x}\right].
\end{dcases}
\end{equation}
The hypothesis $w=v$ is justified considering the limit as the parameter $c$ in the third equation of the complete model~\eqref{eqn:system} approaches the infinity and, as proof of this statement, in~\cite{Moschetta 2019} is shown that, by increasing the $c$ value and taking into account the corresponding wave speed approximation for the cancerous cells density achieved in the full model case, a convergence towards the approximated wave speed provided by the simplified model~\eqref{eqn:system2} is appreciable.

It is also important to notice that the assumption $w=v$ produces a quantitative mismatch concerning the wave speed exhibited by the full model with respect to the reduction, explainable detecting that the reduced model cannot lean on two diffusion mechanisms as it happens for the full system. As a result, the propagating fronts speed computed for the reduction is smaller. On the other hand, talking about the qualitative aspect, the simplified model accomplishes the purpose of correctly reproducing both the heterogeneous and homogeneous configurations, although the \emph{gap} formation is no longer observable: indeed, in order to detect this phenomenon, is mandatory exploiting an independent evolution for the lactic acid concentration.

\subsection{The numerical algorithm}
As regards the numerical strategy, it is back to what has been explained in~\cite{Moschetta 2019}. For the sake of convenience, we provide the reader with the highlights underlying the system discretization. 

\begin{figure}[!ht]
\includegraphics[width=.89\textwidth]{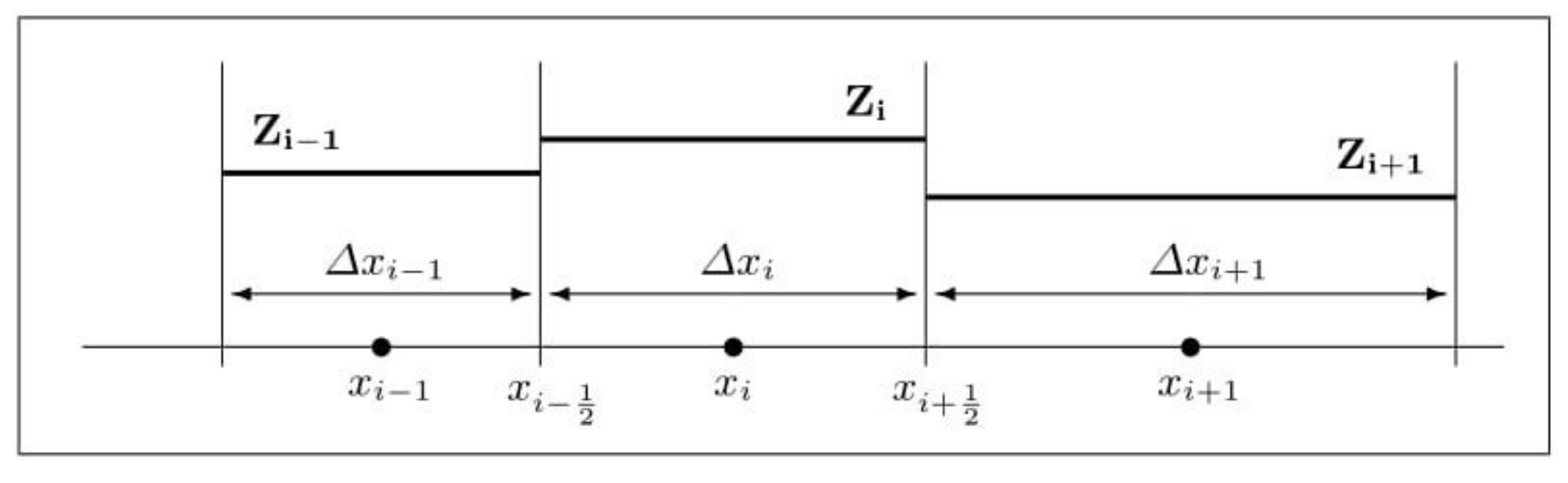}
\caption{Piece-wise constant reconstruction on nonuniform mesh.}
\label{fig:spatialmesh}
\end{figure}

We employ a cell-centered finite volume approximations for the spatial discretization (see~\cite{Wesseling 2001}, for example) and proceed by considering a nonuniform grid. Thus, let $Z_{i}=[x_{i-\frac12},x_{i+\frac12})$ be the finite volume centered at $x_i=\dfrac{x_{i-\frac12}+x_{i+\frac12}}{2}$, for $i=1,2,..., N$, where $N$ is a fixed number of vertices on the one-dimensional grid.\\ Let us assume that $\Delta x_{i}=|x_{i+\frac12}-x_{i-\frac12}|$ is the spatial grid size, from which $|x_i-x_{i-1}|=\dfrac{\Delta x_{i-1}}2+\dfrac{\Delta x_i}2$ is the length for an interfacial interval (see Figure~\ref{fig:spatialmesh}).

The finite volume integral version for the healthy cells density in~\eqref{eqn:system2} leads to
\begin{equation*}
\frac1{\Delta x_i} \!\int_{Z_i} \frac{\partial u}{\partial t}(x,t)\,dx = \frac1{\Delta x_i} \!\int_{Z_i} \!u(x,t)\bigl(1-u(x,t)\bigr)\,dx - \frac{d}{\Delta x_i} \!\int_{Z_i} \!u(x,t) v(x,t)\,dx\,,
\end{equation*}
that, by exploiting the standard notation $\displaystyle{u_i(t) \simeq \frac1{\Delta x_i} \!\int_{Z_i} \!u(x,t)\,dx}$, becomes
\begin{equation}
\label{eqn:finitevolume1}
\frac{d}{dt} u_i(t) = u_i(t) \bigl(1-u_i(t)\bigr) - d\,u_i(t) v_i(t)\,.
\end{equation}
The equation for the tumour cells density in~\eqref{eqn:system2} reads as
\begin{equation*}
\begin{split}
\frac1{\Delta x_i} \!\int_{Z_i} \frac{\partial v}{\partial t}(x,t)\,dx = & \,\frac{r}{\Delta x_i} \!\int_{Z_i} \!v(x,t)\bigl(1-v(x,t)\bigr)\,dx \\
&+ \frac{D}{\Delta x_i} \!\int_{Z_i} \frac{\partial }{\partial x}\left[\bigl(1-u(x,t)\bigr) \frac{\partial v}{\partial x}(x,t)\right] dx
\end{split}
\end{equation*}
where the finite volume integral average of the diffusion term is to be properly approached, as follows, 
\begin{equation}
\label{eqn:finitevolume2}
\begin{split}
& \frac{D}{\Delta x_i} \left[ \bigl(1-u(x_{i+\frac12},t)\bigr) \frac{\partial v}{\partial x}(x_{i+\frac12},t) - \bigl(1-u(x_{i-\frac12},t)\bigr) \frac{\partial v}{\partial x}(x_{i-\frac12},t) \right]\\
&\simeq \frac{D}{\Delta x_i} \left[ \frac{\bigl(1-u_i(t)\bigr)\Delta x_i+\bigl(1-u_{i+1}(t)\bigr)\Delta x_{i+1}}{\Delta x_i+\Delta x_{i+1}} \cdot \frac{v_{i+1}(t)-v_i(t)}{\dfrac{\Delta x_i}2+\dfrac{\Delta x_{i+1}}2} \right.\\
& \qquad\qquad \left. - \,\frac{\bigl(1-u_{i-1}(t)\bigr)\Delta x_{i-1}+\bigl(1-u_i(t)\bigr)\Delta x_i}{\Delta x_{i-1}+\Delta x_i} \cdot \frac{v_i(t)-v_{i-1}(t)}{\dfrac{\Delta x_{i-1}}2+\dfrac{\Delta x_i}2} \right],
\end{split}
\end{equation}
where the approximations for the interfacial quantities are realized by means of weighted averages whose weights are the size of the adjacent finite volumes, thus $\Delta x_i / \Delta x_{i+1}$ and $\Delta x_{i-1} / \Delta x_i$ are employed at the interfaces $x_{i+\frac12}$ and $x_{i-\frac12}\,$, respectively. The first order derivatives of $v(x,t)$ are discretized through an upwind formula which relies on the function evaluations at the neighboring vertices.

From now on, we simply impose that the quantity $\Delta x_i$ is constant, so that $\Delta x_i=\Delta x$ for all $i=1,2,...,N$. That is why, from~\eqref{eqn:finitevolume2} the semi-discrete version for the equation of cancerous cells density reads as
\begin{equation*}
\begin{split}
\frac{d}{dt} v_i(t) = r\,v_i(t)\bigl(1-v_i(t)\bigr) + \frac{D}{\Delta x} & \left[ \frac{\bigl(1-u_i(t)\bigr)\!+\!\bigl(1-u_{i+1}(t)\bigr)}{2} \cdot \frac{v_{i+1}(t)-v_i(t)}{\Delta x} \right.\\
& \;\left. - \,\frac{\bigl(1-u_{i-1}(t)\bigr)\!+\!\bigl(1-u_i(t)\bigr)}{2} \cdot \frac{v_i(t)-v_{i-1}(t)}{\Delta x} \right]
\end{split}
\end{equation*}
which can be rearranged to get
\begin{equation}
\label{eqn:finitevolume2bis}
\begin{split}
\frac{d}{dt} v_i(t) = r\,v_i(t)\bigl(1-v_i(t)\bigr) + \frac{D}{\Delta x^2} & \biggl[ \bigl(1-u_i(t)\bigr)\bigl(v_{i+1}(t)-2\,v_i(t)+v_{i-1}(t)\bigr)\\
&\; - \frac12 \bigl(v_{i+1}(t)-v_i(t)\bigr) \bigl(u_{i+1}(t)-u_i(t)\bigr)\\
&\; - \frac12 \bigl(v_i(t)-v_{i-1}(t)\bigr) \bigl(u_i(t)-u_{i-1}(t)\bigr) \biggr].
\end{split}
\end{equation}
As already pointed out in~\cite{Moschetta 2019}, we stress that the approximation~\eqref{eqn:finitevolume2bis} produces a discrete Laplace operator and extra terms consisting of products of upwind discretizations, arising from the degenerate diffusion in the second equation of~\eqref{eqn:system2}. It is important to notice that the finite volume strategy allows to split the diffusion by autonomously choosing the first and second order contributions, while, as in the case of finite difference schemes~\cite{Quarteroni 2014}, a central discretization for the first order terms would be required, so causing a less stable scheme.

Finally, for the time discretization of ~\eqref{eqn:finitevolume1} and \eqref{eqn:finitevolume2bis}, we adopt a semi-implicit strategy considering a fixed time step $\Delta t\,$, thus $\Delta t=|t^{n+1}-t^{n}|$, for $n=0,1,...$. The reaction terms are treated explicitly, while the differential terms on the right-hand sides are approximated implicitly, as follows,
\begin{equation}
\label{eqn:timescheme}
\begin{dcases}
u^{n+1}_i = & u^n_ i + \Delta t \biggl[ u^n_i \bigl(1-u^n_i \bigr) - d\,u^n_i v^n_i \biggr]\\
v^{n+1}_i = & v^n_ i + r \Delta t\,v^n_i \bigl(1-v^n_i \bigr)\\
& +\,D \frac{\Delta t}{\Delta x^{2}} \biggl[ \bigl(1-u^{n+1}_i \bigr) \bigl(v^{n+1}_{i+1}-2\,v^{n+1}_i+v^{n+1}_{i-1} \bigr)\\
&\qquad\qquad\; -\,\frac12 \bigl(v^{n+1}_{i+1}-v^{n+1}_i \bigr) \bigl(u^{n+1}_{i+1}-u^{n+1}_i \bigr)\\
&\qquad\qquad\; -\,\frac12 \bigl(v^{n+1}_i-v^{n+1}_{i-1} \bigr) \bigl(u^{n+1}_i-u^{n+1}_{i-1} \bigr) \biggr]
\end{dcases}
\end{equation}
and Neumann-type boundary conditions $u_1^n=u_2^n$ and $v_1^n=v_2^n$, for $n=1,2,...$ are implemented. 

\subsection{Simulations results}
In continuity with~\cite{Moschetta 2019}, we perform simulations aimed at better characterizing the solutions produced by~\eqref{eqn:system2}. First of all, we want to figure out if the corresponding traveling waves exhibit a \textit{sharp-type} or \textit{front-type} trend. Technically, taking as main guideline the traveling waves problem defined by the following one-dimensional, degenerate, reaction-diffusion equation, 
\begin{equation}
\label{eqn:degeq}
\frac{\partial v}{\partial t} = \frac{\partial}{\partial x} \left[F(v) \frac{\partial v}{\partial x}\right]+g(v) \qquad \mbox{with} \quad (x,t) \in (\mathbb{R} \times \mathbb{R}^{+}),
\end{equation}
where $g(v)$ is a Fisher-KPP type reaction term, $F(v)$ is the degenerate diffusion such that $v \in [0,1]$ and $F^{\prime}(0)\not= 0$, then, the definition of \textit{sharpness}, according to~\cite{Sanchez 1995}, reads as 
\begin{definition}[sharp-type front]
\label{def:sharpdef}
If there exist a value of the wave speed $s$, let us call it $s^{\ast}$, and a value of $\xi$, let it be $\xi^{\ast} \in (-\infty,+\infty]$, such that $\phi(x-s^{\ast}t)=\phi(\xi)$, satisfying
\begin{enumerate}
\item $F(\phi)\phi^{\prime\prime}+s^{\ast}\phi^{\prime}+ F^{\prime}(\phi)[\phi^{\prime}]^{2}+g(\phi)=0 
\quad \forall \xi \in (-\infty,\xi^{\ast})$,
\item $\phi(-\infty)=1, \quad \phi(\xi^{\ast})=0 \quad \mbox{and} \quad \phi^{\prime}<0 \quad \forall \xi \in (-\infty,\xi^{\ast})$,
\item $\phi^{\prime}(\xi^{\ast})=-s^{\ast}/F^{\prime}(0) \quad \mbox{and} \quad \phi(\xi)=0 \quad \forall \xi \in (\xi^{\ast},+\infty]$,
\end{enumerate}
where the superscript is meant to denote differentiation with respect to $\xi$, then the function $v(x,t)=\phi(x-s^{\ast}t)$ is called a traveling wave solution of \textit{sharp-type} for~\eqref{eqn:degeq}. 
\end{definition} 
We point out that the other possibility allowed, happens when the function $v(x,t)$ turns out to be a traveling wave of \textit{front-type}, whose typical smoother trend makes this front to be known as a \textit{smooth-type} wave as well. The former statement about the smoothness of the front-type traveling waves, is easily understandable thinking about the implications framed by the Definition~\ref{def:sharpdef}. As a consequence, indeed, a sharp-type wave attains the equilibrium located in $\mathbf{E}=0$ in a finite time $\xi^{\ast}$, with negative slope $\phi^{\prime}(\xi^{\ast})=-s^{\ast}/F^{\prime}(0)$~\cite{Malaguti 2003}, thus resulting a discontinuous derivative in $\xi^{\ast}$, since the left derivative tends to $\phi^{\prime}(\xi^{\ast-}) \not= 0$, while the right derivative tends to $\phi^{\prime}(\xi^{\ast+})= 0$~\cite{Sanchez 1995}. By contrast, a smooth-type front exhibits a continuous derivative in $\xi^{\ast}$. This last observation provides us with a useful tool in order to quickly, qualitatively detect the distinctive trend for a given traveling wave, especially when the dynamics is ruled by more complex configurations with respect to~\eqref{eqn:degeq}, as it happens for a system of equations. As concerns the scalar case, in which the problem~\eqref{eqn:degeq} is framed, theoretical results~\cite{Malaguti 2003, Sanchez 1994 first, Sanchez 1995, Sanchez 1996} are available for ensuring the existence and uniqueness of sharp/smooth-type traveling waves, provided that some hypotheses about the regularity of the $v$-dependent functions $F$ and $g$ are satisfied; other results are achieved in~\cite{Sanchez 1994 second} for a specific choice of $F$ and in~\cite{Sanchez 1997} if $g$ is a generalization of the Nagumo equation. 

As far as the strictly theoretical framework about~\eqref{eqn:degeq}, we make some considerations in the next section, when the one-equation reduction for the Gatenby-Gawlinski model is introduced; right now, instead, we focus on the sharpness for~\eqref{eqn:system2}, without neglecting observations concerning~\eqref{eqn:system} too. 

\begin{figure}[!ht]
\includegraphics[width=.55\textwidth]{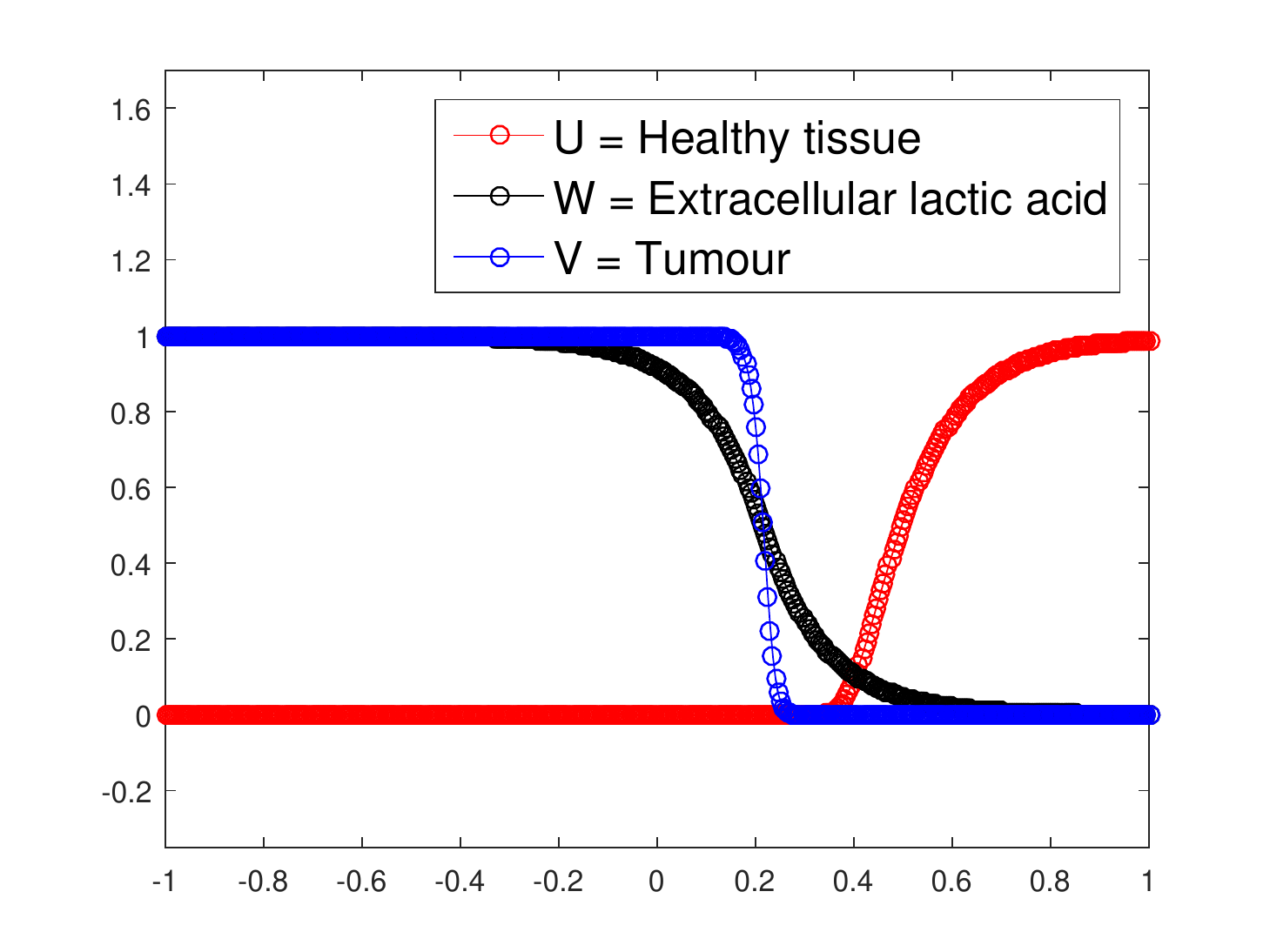}
\caption{Simulations results in the full model case \eqref{eqn:system} for the homogeneous invasion. The parameters used are listed in Table~\ref{tab:parameters}.}
\label{fig:3eqprofiles}
\end{figure}

The first step, indeed, consists in evaluating the traveling fronts arising from the full model: in this regard, we take as a sample (see Figure~\ref{fig:3eqprofiles}) the results related to the homogeneous invasion considered in~\cite{Moschetta 2019}, but initialized with the Riemann problem whose states are suitable stationary points~\cite{McGillen 2014} of the full model. The parameters used for the experiment are listed in Table~\ref{tab:parameters}; moreover, $T$ is the final time instant, while the spatio-temporal mesh is realized by fixing $\Delta x = 0.005$ and $\Delta t = 0.005$. Talking about the numerical algorithm, we exploit the strategy previously described provided with the equation for the lactic acid concentration. 

\begin{table}[!ht]
\caption{Numerical default values for the parameters in the complete model case.}
\label{tab:parameters}
\begin{tabular}{*{4}{c}}
\toprule
$\mathbf{d}$ & $\mathbf{r}$ & $\mathbf{D}$  & $\mathbf{T}$\\
\midrule
$12.5$ & $1$ & $4 \cdot 10^{-5}$  & $20$\\
\bottomrule
\end{tabular}
\end{table}

\begin{figure}[!ht]
\subfloat[][\emph{healthy cells fronts}]
{\includegraphics[width=.49\columnwidth]{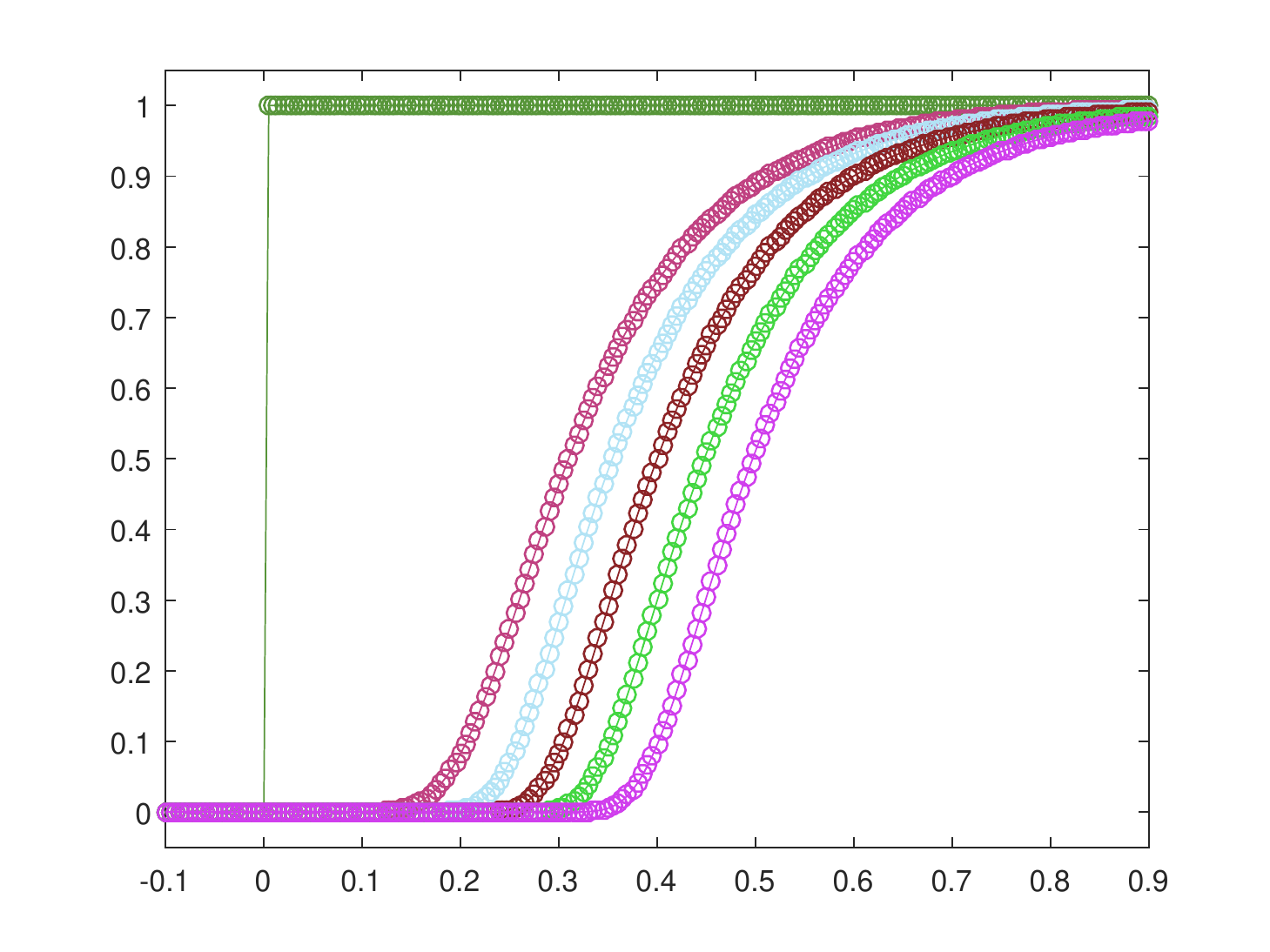}}
\subfloat[][\emph{tumour cells fronts}]
{\includegraphics[width=.49\columnwidth]{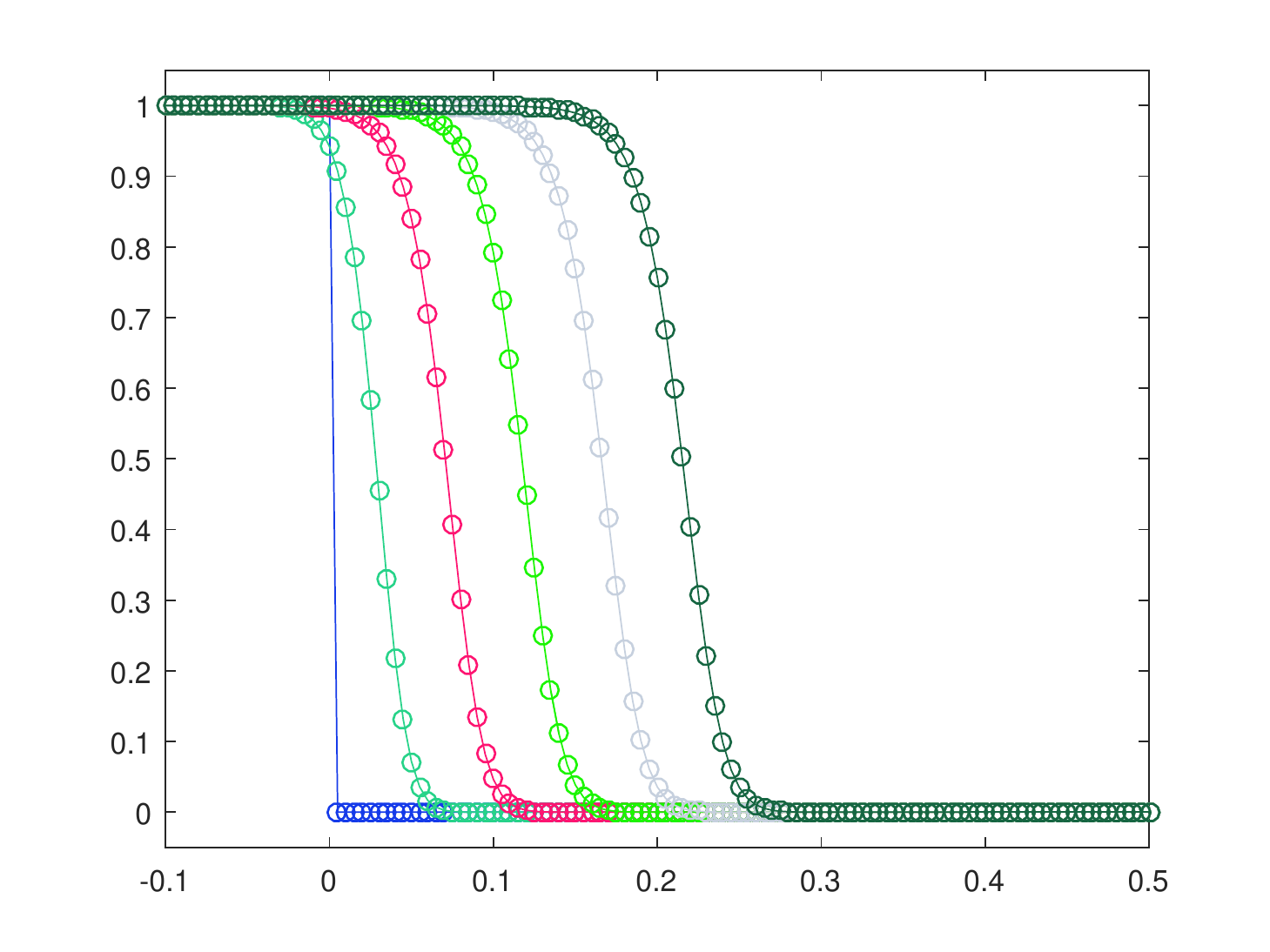}}
\caption{Front evolution zoom-in for the healthy cells density (A) and the tumour cells density (B).}
\label{fig:3eqsmooth}
\end{figure}

Now, in order to get information about the shape of the fronts, we realize a zoom-in for both the healthy cells density, shown in Figure~\ref{fig:3eqsmooth}(A), and the tumour cells density, available in Figure~\ref{fig:3eqsmooth}(B), plotted at equally spaced time instants. This qualitative analysis clearly proves the traveling waves to be smooth-type for the complete Gatenby-Gawlinski model. Concerning the lactic acid concentration, we have omitted to report data, due to the similarity with the evolution of cancerous cells density.

\begin{figure}[!ht]
\subfloat[][\emph{heterogeneous invasion}]
{\includegraphics[width=.49\columnwidth]{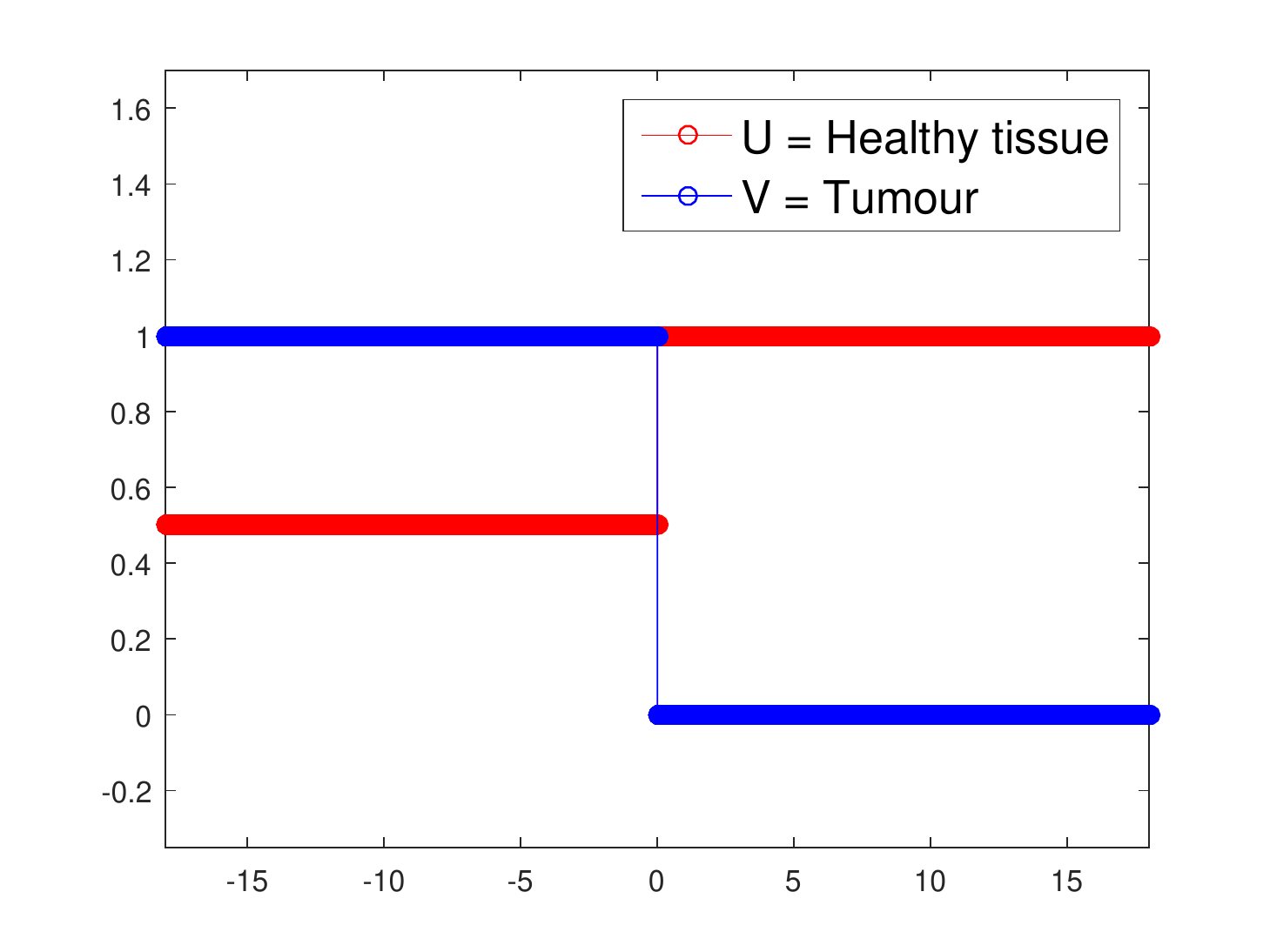}}
\subfloat[][\emph{homogeneous invasion}]
{\includegraphics[width=.49\columnwidth]{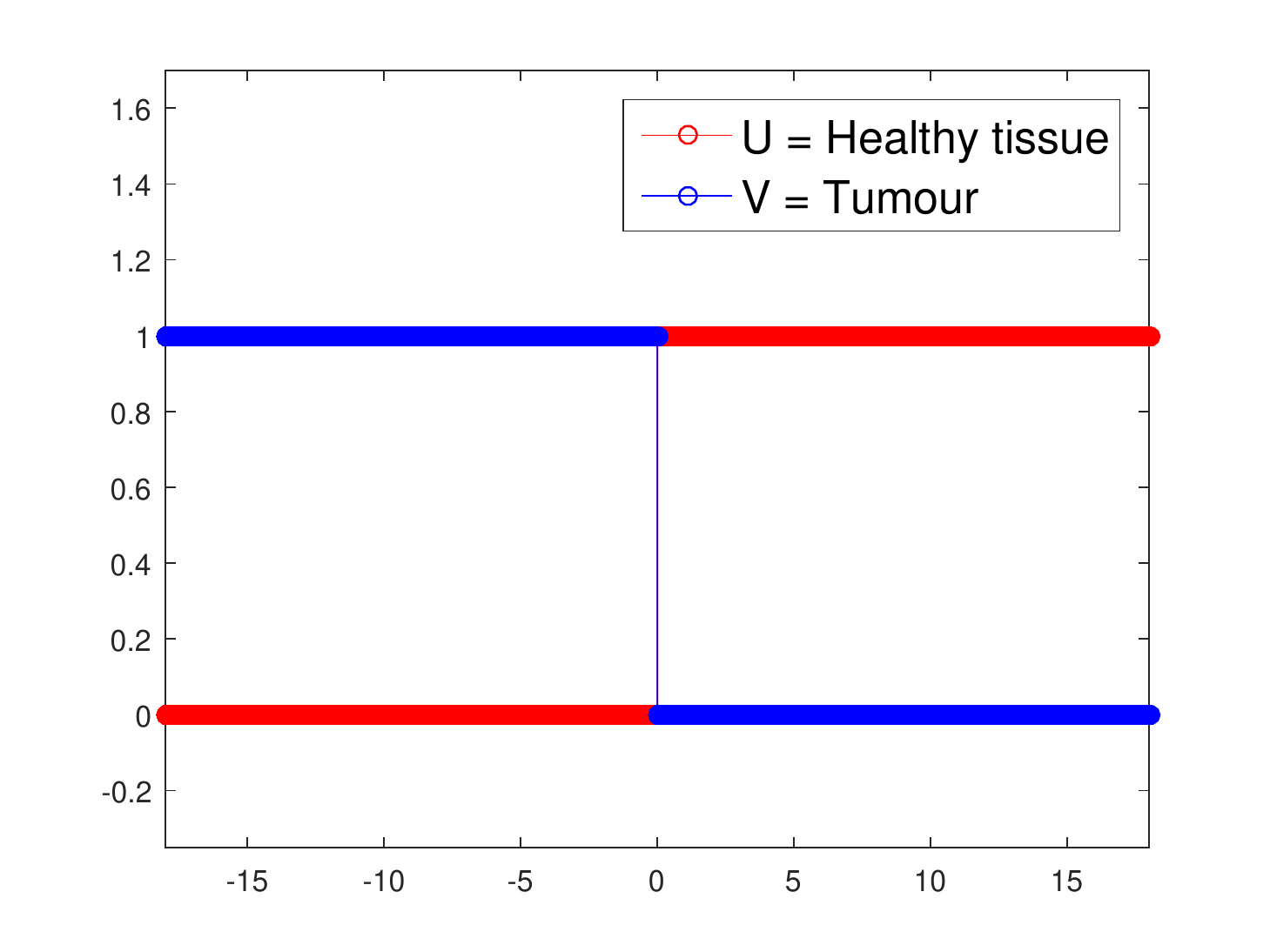}}
\caption{Riemann problem for the heterogeneous case (A) and the homogeneous case (B). The parameters used are listed in Table~\ref{tab:parametersReduc}.}
\label{fig:Riemann}
\end{figure}

Let us proceed by focusing on the two-equations-based model~\eqref{eqn:system2}. Before going ahead with the numerical simulations, it turns out useful making a further simplification, allowing to normalize the coefficient $D$ for the second equation (it should be noted that the resulting spatial window is wider). This goal is accomplished imposing the following rescaling
\begin{equation}
\label{eqn:rescaling}
\sqrt D \frac{\partial }{\partial x} = \frac{\partial }{\partial y},
\end{equation}
under which, by renaming the variable $y$ to $x$, it is possible to get a two-parameters-dependent reduction that reads as 
\begin{equation}
\label{eqn:system3}
\begin{dcases}
\frac{\partial u}{\partial t} = u (1-u) - d u v\\
\frac{\partial v}{\partial t} = r v (1-v) + \frac{\partial}{\partial x} \left[(1-u) \frac{\partial v}{\partial x}\right].
\end{dcases}
\end{equation}
Afterwards, it is helpful to compute the stationary points for~\eqref{eqn:system3}, namely~\eqref{eqn:system2}, and check the related stability, so that the final outcome looks like
\begin{itemize}
\item $\mathbf{E_{0}}=(0,0)$, absence of species, unstable; 
\item $\mathbf{E_{1}}=(1,0)$, healthy state, unstable; 
\item $\mathbf{E_{2}}=(0,1)$, homogeneous state, stable if $d>1$ and unstable if $d<1$; 
\item $\mathbf{E_{3}}=(1-d,1)$, heterogeneous state, stable if $d<1$ and unstable if $d>1$. 
\end{itemize}
By taking advantage of this report, we consider the Riemann problem whose states are $\mathbf{E}=\bigl((1-d)^{+},1\bigr)$ at the left and $\mathbf{E_{1}}=(1,0)$ at the right, as guideline for the initial profiles to be selected for performing numerical experiments (see Figure~\ref{fig:Riemann}(A) for the heterogeneous invasion, with $d=0.5$, and Figure~\ref{fig:Riemann}(B) for the homogeneous invasion, with $d=2$). 

\begin{figure}[!ht]
\subfloat[][\emph{heterogeneous invasion}]
{\includegraphics[width=.49\columnwidth]{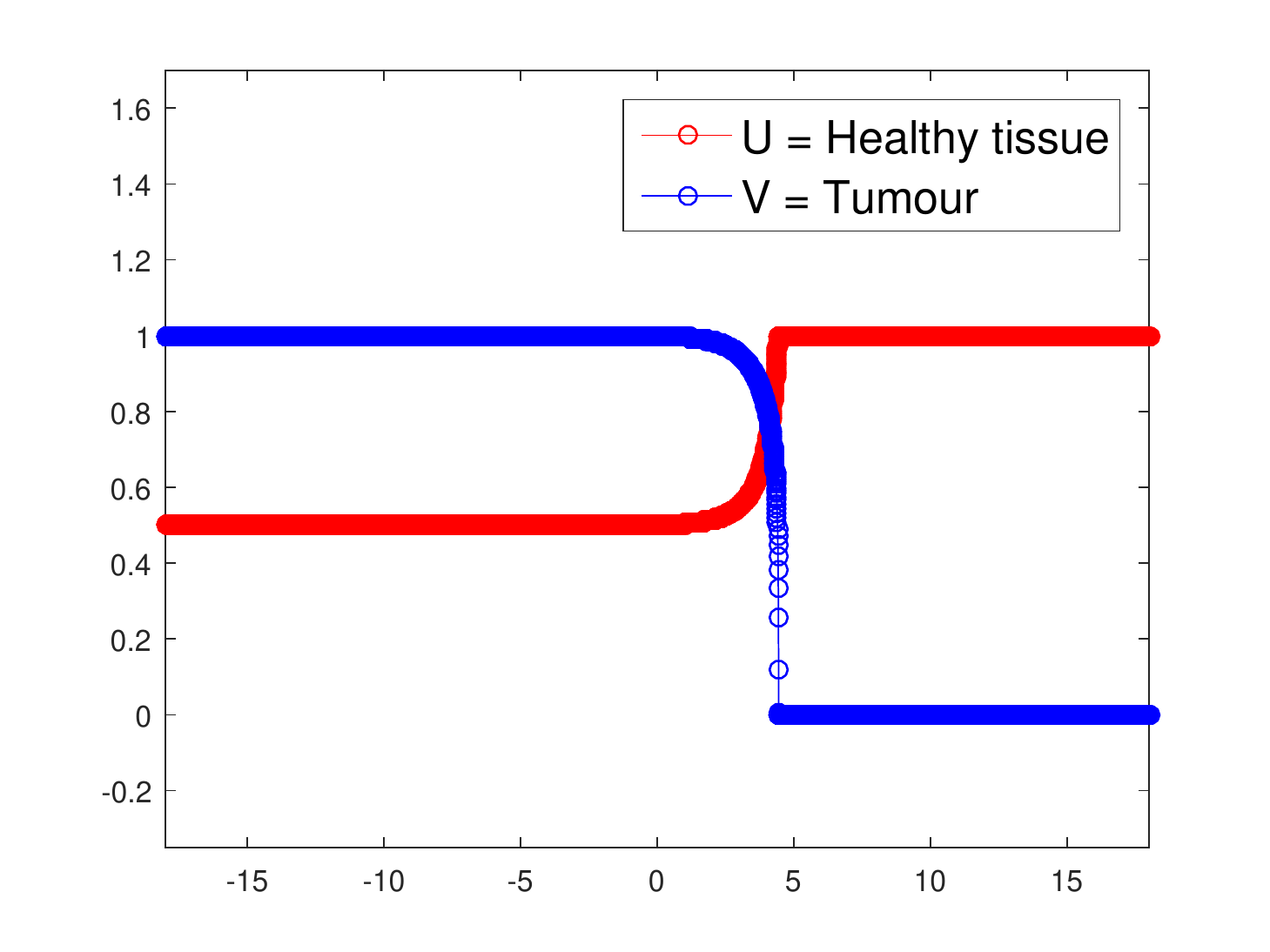}}
\subfloat[][\emph{homogeneous invasion}]
{\includegraphics[width=.49\columnwidth]{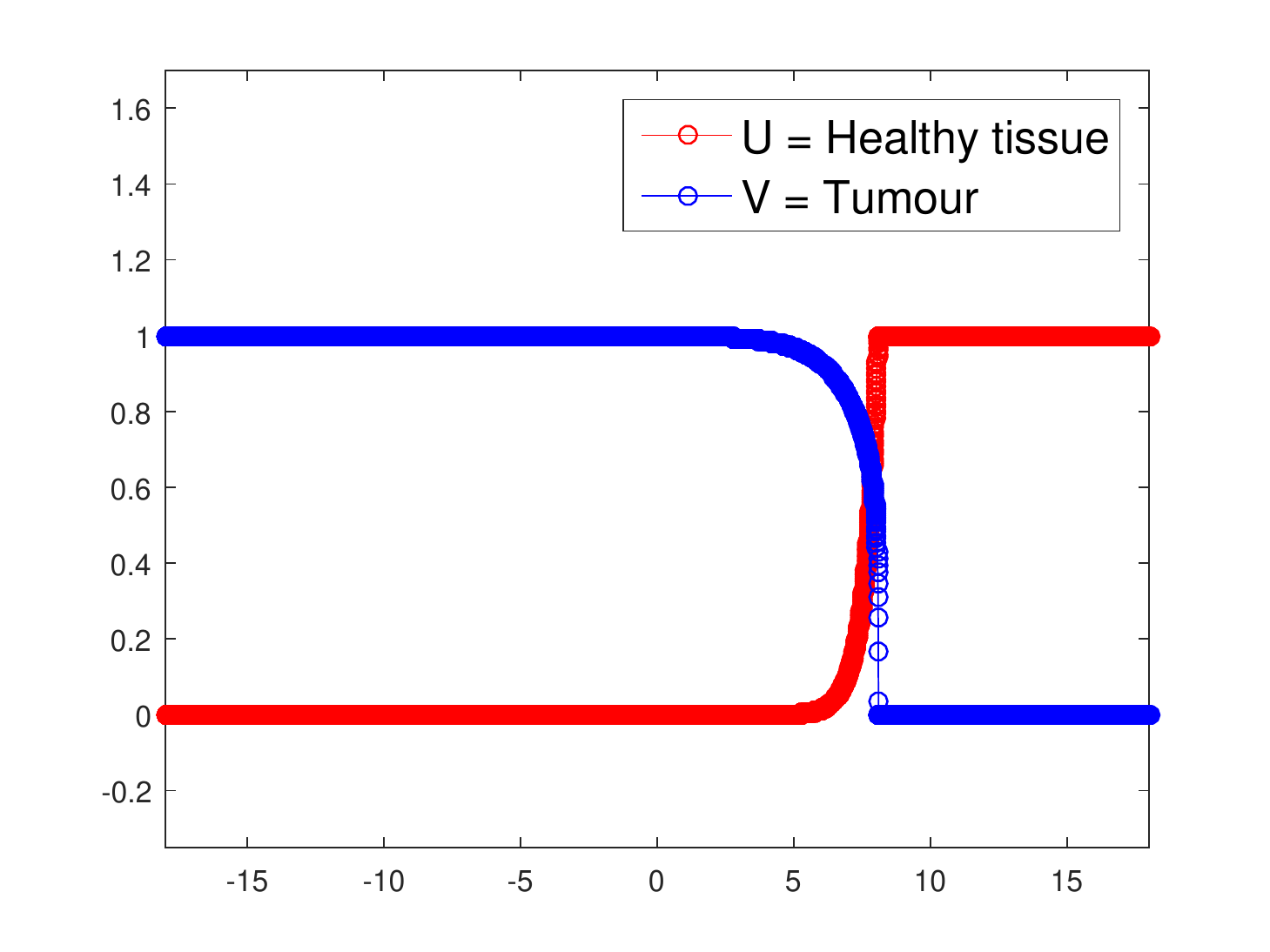}}
\caption{Numerical solutions for the heterogeneous case (A) and the homogeneous case (B). The parameters used are listed in Table~\ref{tab:parametersReduc}.}
\label{fig:solutions}
\end{figure}

It is important to notice that the structure of the state $\mathbf{E}$ is a consequence of the $d$-dependent stability for the equilibria $\mathbf{E_{2}}$ and $\mathbf{E_{3}}$, so that the resulting state $\mathbf{E}$ proves itself to be stable for every choice of eligible $d$. Being $\mathbf{E_{1}}$ unstable, the propagating front arising from the resulting Riemann data, travels towards the right-hand side. In Table~\ref{tab:parametersReduc}, all the parameters employed to carry out numerical simulations are shown. 
\begin{table}[!ht]
\caption{Numerical default values for the parameters involved in the simplified model.}
\label{tab:parametersReduc}
\begin{tabular}{*{5}{c}}
\toprule
$\mathbf{d}$ & $\mathbf{r}$ & $\mathbf{\Delta x}$ & $\mathbf{\Delta t}$  & $\mathbf{T}$\\
\midrule
$\{0.5,2\}$ & $1$ & $0.005$ & $0.005$ & $20$\\
\bottomrule
\end{tabular}
\end{table} 
The results produced by numerically investigating the system~\eqref{eqn:system3} are depicted in Figure~\ref{fig:solutions}(A) for the heterogeneous invasion and Figure~\ref{fig:solutions}(B) for the homogeneous invasion: we recognize that the two-parameters reduction correctly catches the trends expected in both the cases under analysis. By making a comparison with the analogous plots obtained in~\cite{Moschetta 2019}, we point out that we have currently employed a smaller $\Delta t$ and defined a wider spatial window to frame the front evolution, since the waves are traveling faster as a consequence of the rescaling~\eqref{eqn:rescaling}. These two changes explain the small differences concerning the steepness and the graphical display of the fronts, being the ones reported in Figure~\ref{fig:solutions} less steep and better graphically depicted. Regardless of these points, the typical trends characterizing  cancer invasions in the Gatenby-Gawlinski model, are definitely, qualitatively preserved by simulations results.

\begin{figure}[!ht]
\subfloat[][\emph{heterogeneous case}]
{\includegraphics[width=.49\columnwidth]{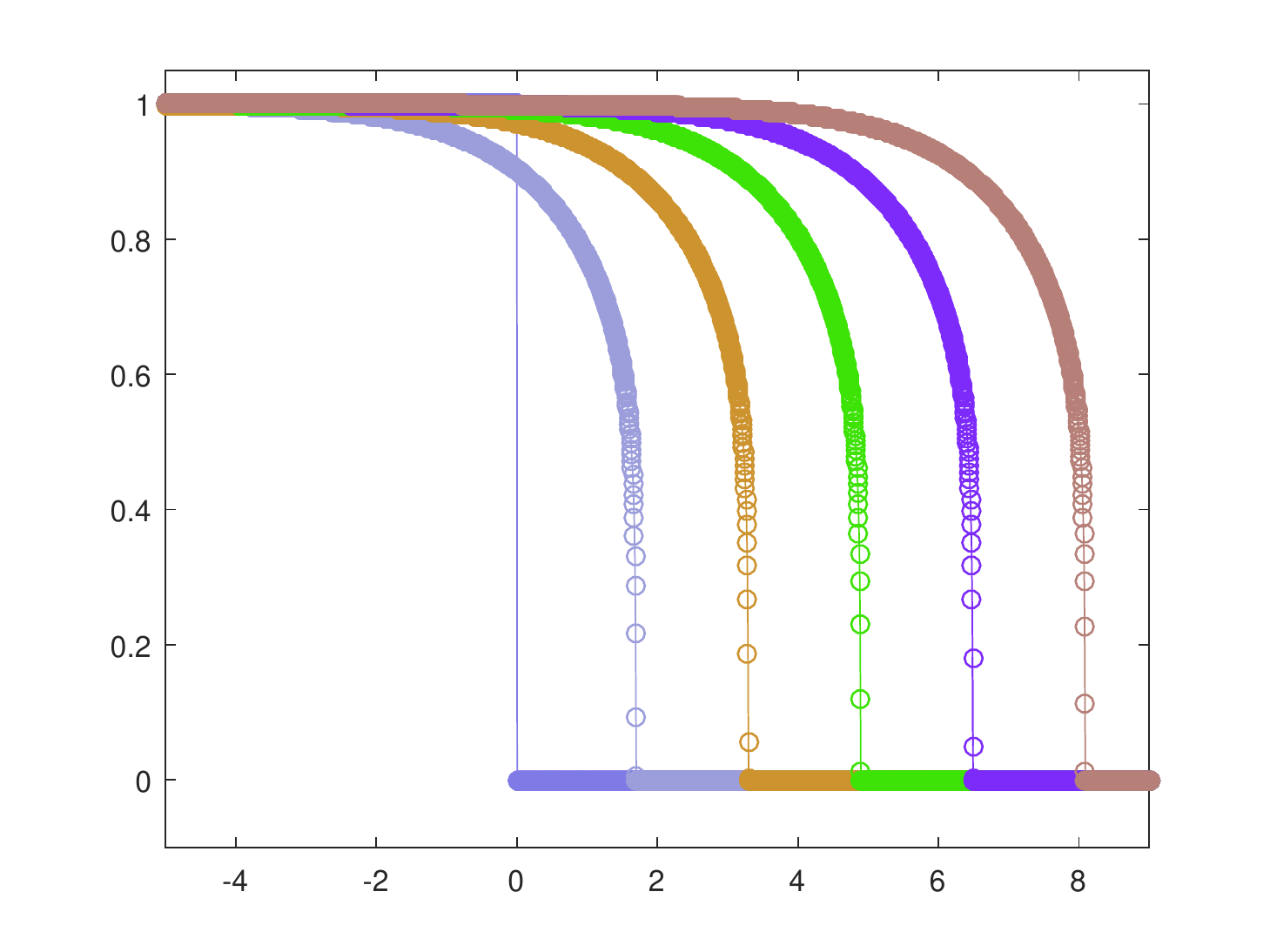}}
\subfloat[][\emph{homogeneous case}]
{\includegraphics[width=.49\columnwidth]{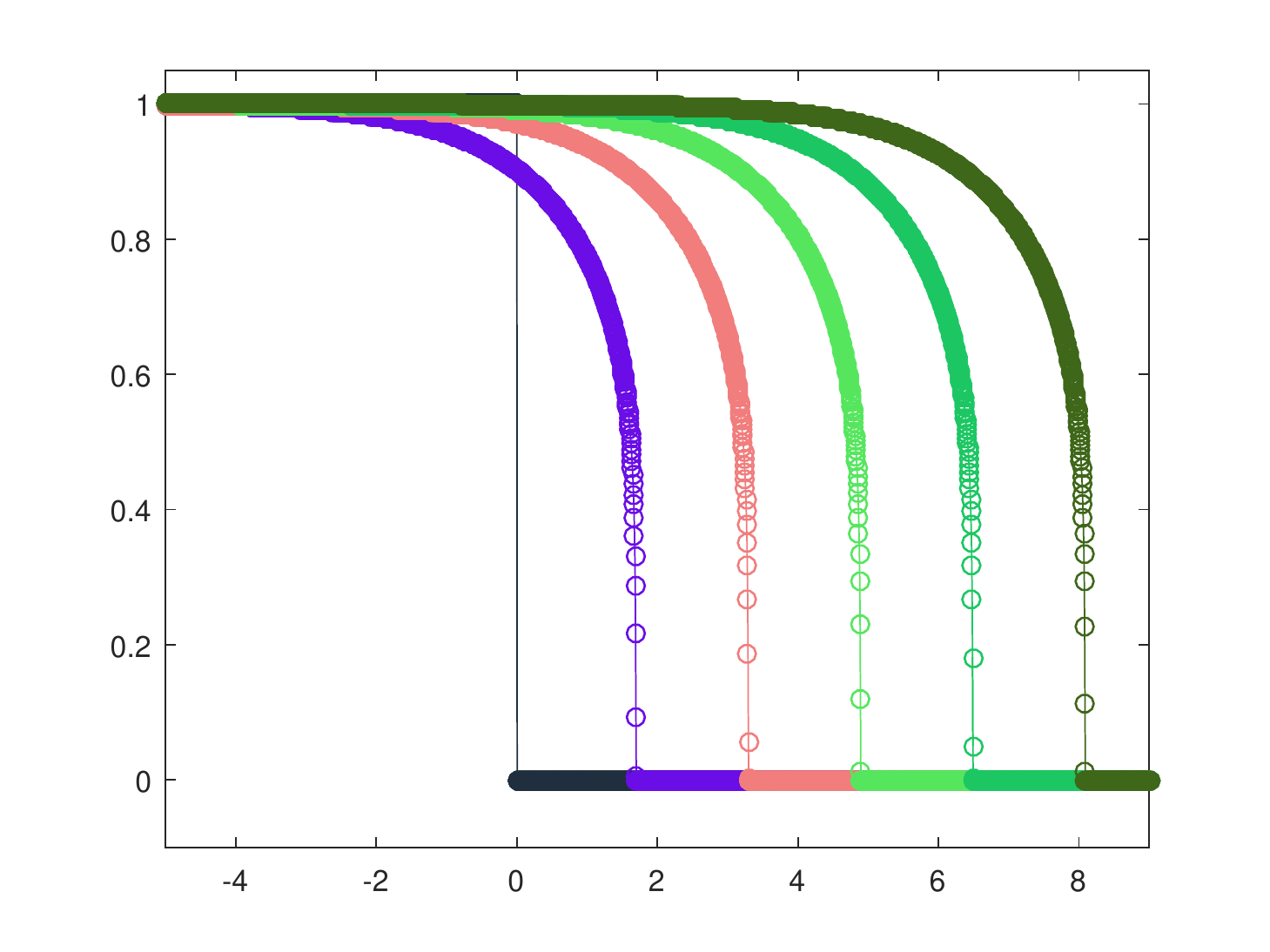}}
\caption{Front evolution zoom-in for the tumour cells density in the heterogeneous case (A) and the homogeneous case (B).}
\label{fig:tumour}
\end{figure}

\begin{figure}[!ht]
\subfloat[][\emph{heterogeneous case}]
{\includegraphics[width=.49\columnwidth]{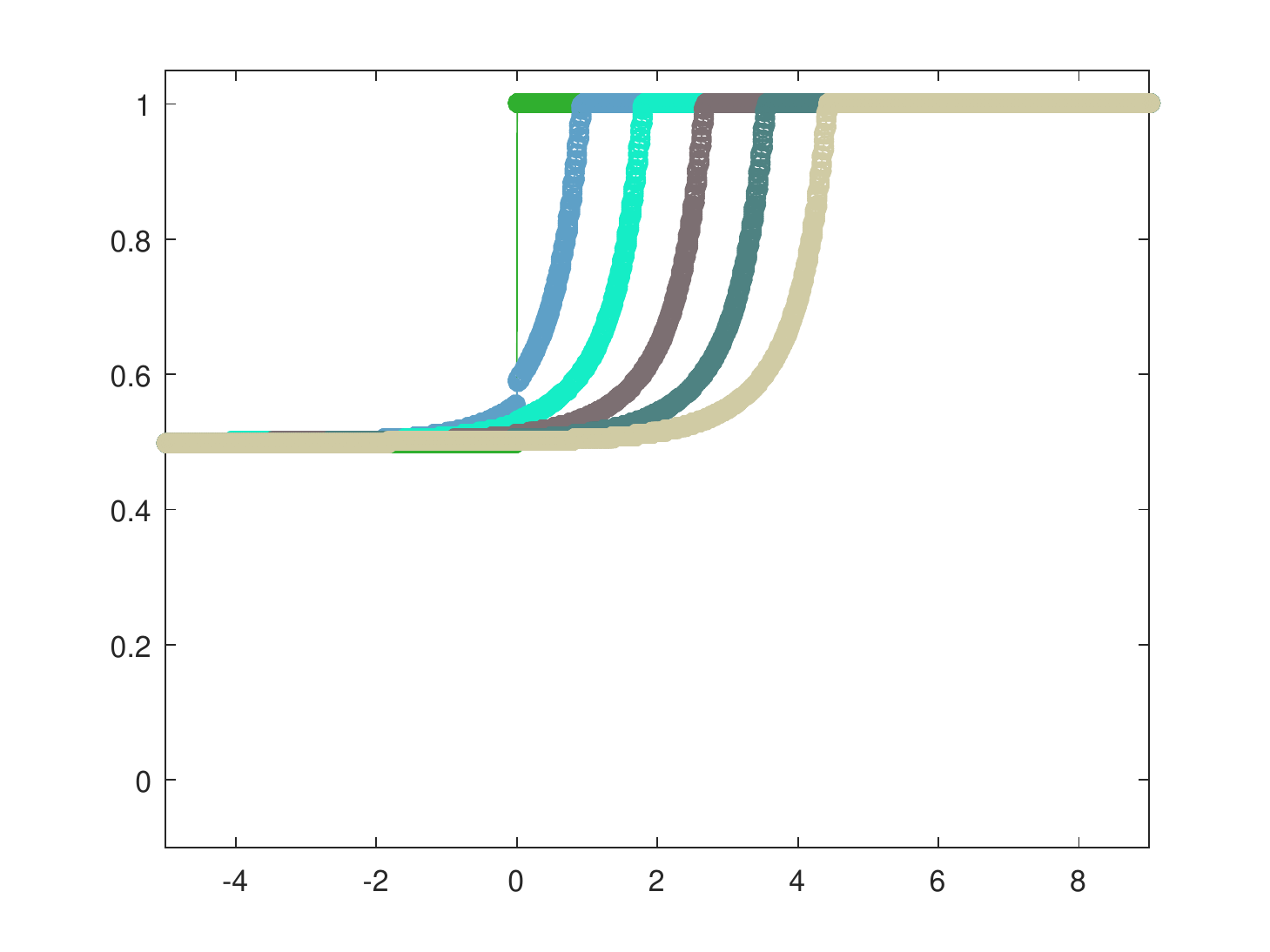}}
\subfloat[][\emph{homogeneous case}]
{\includegraphics[width=.49\columnwidth]{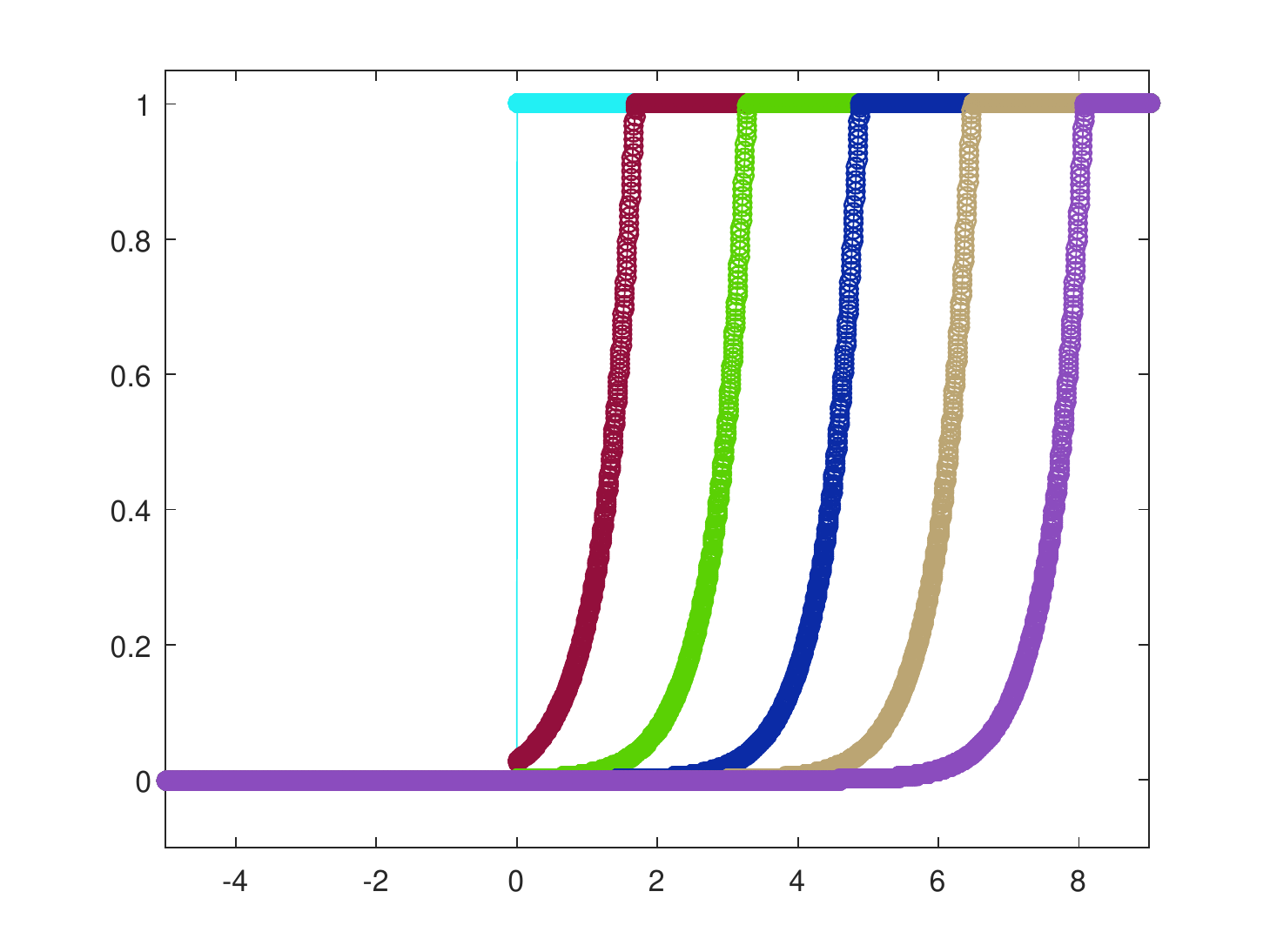}}
\caption{Front evolution zoom-in for the healthy cells density in the heterogeneous case (A) and the homogeneous case (B).}
\label{fig:healthy}
\end{figure}

Next, let us focus on the shape of the traveling waves. Figure~\ref{fig:tumour}(A) and Figure~\ref{fig:tumour}(B) exhibit the front evolution for the tumour cells density in both the invasion configurations, while in Figure~\ref{fig:healthy}(A) and Figure~\ref{fig:healthy}(B) are displayed the corresponding plots for the healthy cells density.
%A check on the plots, clearly ensures the sharpness of the fronts, certified by the detection of a discontinuous derivative:
%this feature is emphasized taking into account the smoothness of the graphs in Figure~\ref{fig:3eqsmooth}.

We conclude this section performing a sensitivity analysis with respect to parameters $r$ and $d$. In order to accomplish this purpose, we consider as unknown the wave speed $s$, whose numerical approximation is made by invoking the space-averaged estimate proposed in~\cite{LeVeque 1990}, already successfully exploited in the Gatenby-Gawlinski model field in~\cite{Moschetta 2019}, to which we refer for the detailed derivation. The final discretized version, providing the approximation for a function $v(x,t)$ over a uniform spatial mesh at time $t^n$ is the following 
\begin{equation}
\label{eqn:numspeed}
s^{n} = \frac{\Delta x}{[\phi]\Delta t} \sum_{i=1}^N \bigl(v_i^n-v_i^{n+1} \bigr),
\end{equation}
where $[\phi]\coloneqq \phi_{+}-\phi_{-}$, being $\phi_{+}$ and $\phi_{-}$ the stationary states of $v(x,t)$. We stress that the estimate~\eqref{eqn:numspeed} is independent from the dynamics of the solutions produced by~\eqref{eqn:system3}.

\begin{figure}[!ht]
\subfloat[][\emph{$r$-parameter sensitivity}]
{\includegraphics[width=.49\columnwidth]{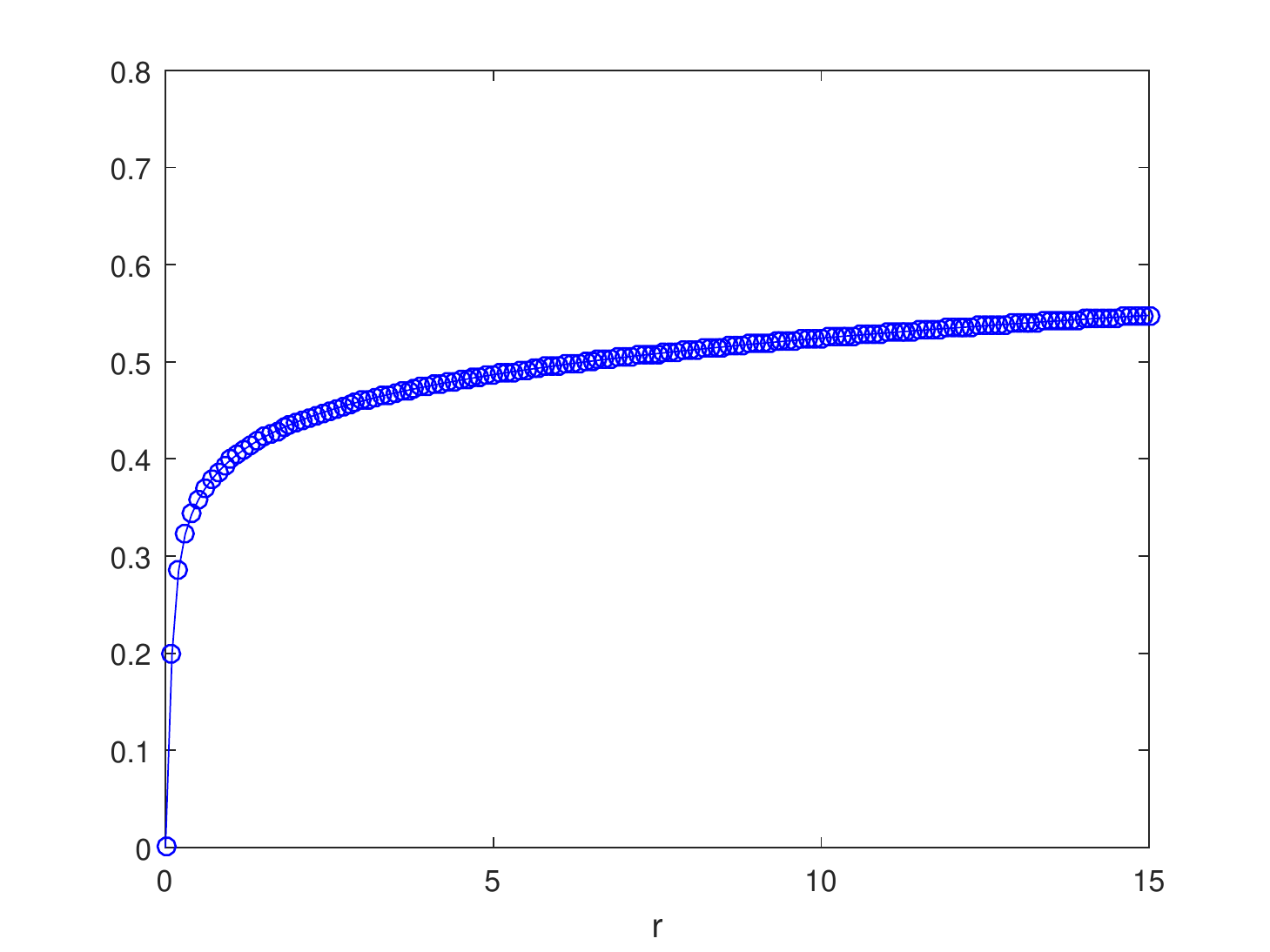}}
\subfloat[][\emph{$d$-parameter sensitivity}]
{\includegraphics[width=.49\columnwidth]{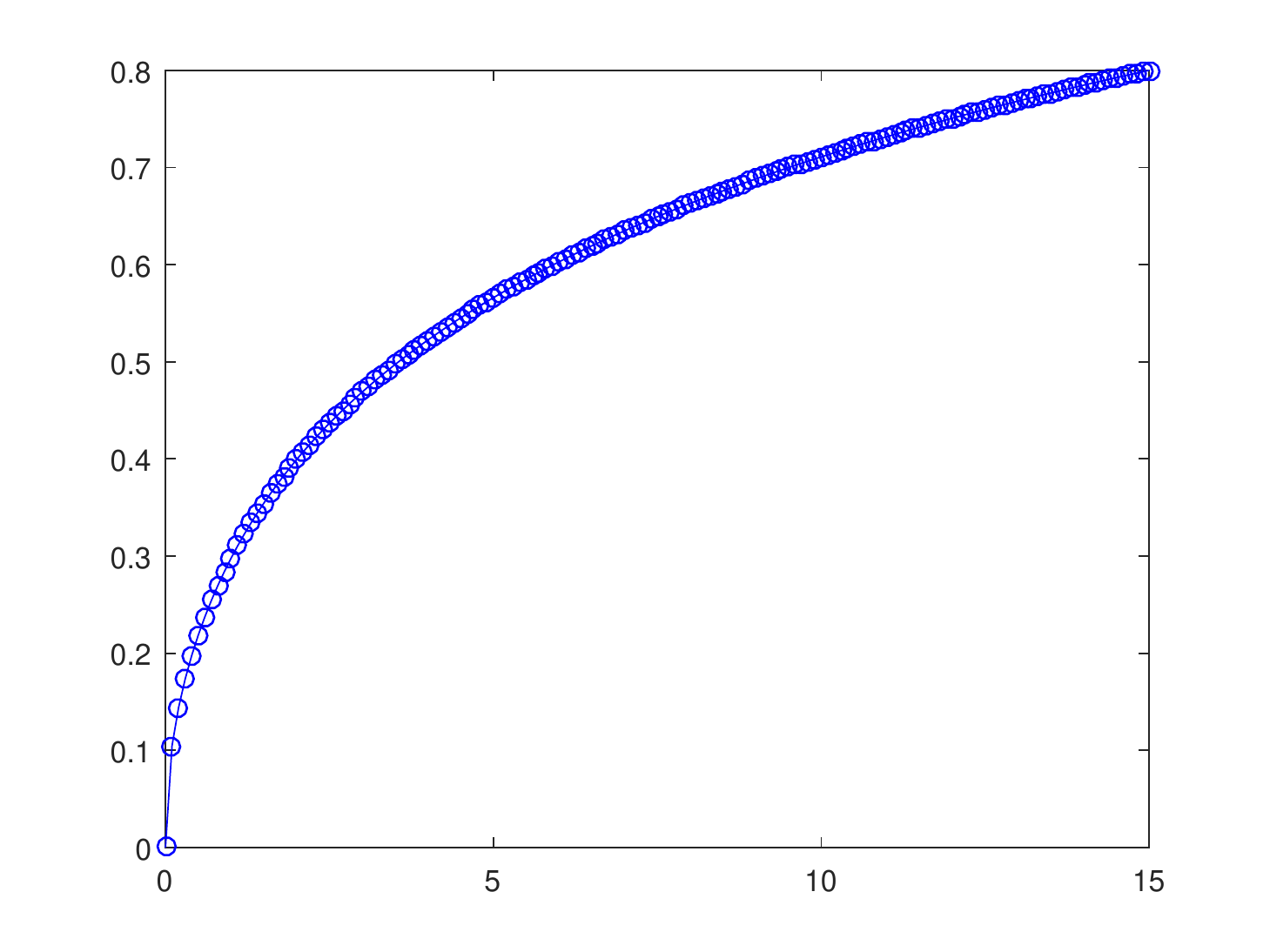}}
\caption{Graph of $s(\cdot,d)$ as a function of the parameter $r \in [0,15]$, with $\Delta r=0.1$, for the homogeneous invasion when $d=2$ (A) and graph of $s(\cdot,r)$ as a function of the parameter $d \in [0,15]$, with $\Delta d=0.1$ and $r=1$ (B).}
\label{fig:speedVSrd}
\end{figure}

In Figure~\ref{fig:speedVSrd}(A), the $r$-dependence is shown taking as a sample the homogeneous invasion (the heterogeneous case exhibits the same qualitative trend). Specifically, for each $r$ value, the corresponding asymptotic wave speed is reported, numerically approximated by means of~\eqref{eqn:numspeed}. The curve so defined is monotone increasing, which is not surprising considering that $r$ is a growth rate for the tumour cells density, so resulting in a faster invasion process. Likewise, focusing on the $d$-dependence, depicted in Figure~\ref{fig:speedVSrd}(B), it follows that cancerous invasion is facilitated as $d$ increases, being this parameter a death rate for the healthy cells due to the interactions with cancerous cells. As a matter of fact, the profile for the wave speed trend is again monotone increasing. % although in the $d$-dependence case a faster growth for the wave speed is recognizable compared to $r$-dependence. This last point shows that cancerous cells, in their invasion, take more advantage from the healthy cells degradation rather than from an increased growth rate for their own reproduction.   

%%%------------------------------------------------------------------------------------------
\section{One-equation-based model reduction}
\label{sec:oneeq}
\subsection{Derivation from the two-equations model}
With the aim of getting a one-equation reduction for the Gatenby-Gawlinski model, we start from the simplified model~\eqref{eqn:system3} and assume the stationarity for the healthy tissue density equation, leading to $u(1-u-d v)=0$; finally we impose that
\begin{equation}
\label{eqn:OneEqHp}
u=(1-dv)^{+}. 
\end{equation}      
As a consequence, the tumour cells equation turns out to show the following structure
\begin{equation}
\label{eqn:OneEqmodelbefore}
\frac{\partial v}{\partial t} = r v (1-v) + \frac{\partial}{\partial x} \left[F(u)\big\rvert_{u=(1-dv)^{+}}\frac{\partial v}{\partial x}\right],
\end{equation} 
where the degenerate diffusion term is a piecewise linear function defined as
\begin{equation}
\label{eqn:DegDiffPLF}
F(u)\big\rvert_{u=(1-dv)^{+}}= 1-u \big\rvert_{u=(1-dv)^{+}}=1-(1-dv)^{+}=\begin{cases} dv & \mbox{if} ~ v \in [0,\frac{1}{d}) \\ 1 & \mbox{if} ~ v \in [\frac{1}{d},1].
\end{cases}
\end{equation}
Our reduction~\eqref{eqn:OneEqmodelbefore} is a degenerate reaction-diffusion equation as in~\eqref{eqn:degeq} and the diffusion $F$ is almost everywhere differentiable, due to the discontinuity located in $v=\frac{1}{d}$ for the derivative. The existence and uniqueness results for traveling waves of sharp/smooth-type available in~\cite{Malaguti 2003}, require at least the pointwise differentiability in $[0,1]$, so that, strictly concerning the theoretical point of view, it is not possible to state the sharpness/smoothness of the fronts for the equation~\eqref{eqn:OneEqmodelbefore}. As a matter of a fact, we keep relying on the numerical assessment in this paper, although the possibility of employing a smooth approximation for by-passing the discontinuous point of $F^{\prime}$, so that enough regularity~\cite{Malaguti 2003, Sanchez 1995} might be ensured, would not seem to jeopardize a theoretical prediction of existence and uniqueness for the fronts in the case of~\eqref{eqn:DegDiffPLF} as well.

We point out that the reduction~\eqref{eqn:OneEqmodelbefore} might be easily rearrenged to become a one-parameter-dependent equation, by a means of the rescaling $\partial/\partial t =r\partial/\partial \tau$. However, taking advantage of the constraint $r=1$ employed for carrying out simulations, it is possible to keep relying on~\eqref{eqn:OneEqmodelbefore} and get a one-parameter dependence anyway, thus leading to the following equation: 
\begin{equation}
\label{eqn:OneEqmodel}
\frac{\partial v}{\partial t} =  v (1-v) + \frac{\partial}{\partial x} \left[F(u)\big\rvert_{u=(1-dv)^{+}}\frac{\partial v}{\partial x}\right].
\end{equation} 

\subsection{The numerical algorithm}   
On the heels of what has already described in Section~\ref{sec:twoeq}, we invoke the same cell-centered finite volume approximation for the spatial discretization of~\eqref{eqn:OneEqmodel} and take care of considering the corresponding versions the piecewise linear diffusion~\eqref{eqn:DegDiffPLF} leads to, specifically we have
\begin{align}
\label{eqn:firstEq}
\frac{\partial v}{\partial t} &=  v (1-v) + d\frac{\partial}{\partial x} \left(v\frac{\partial v}{\partial x}\right) & \mbox{if} ~ v \in \biggl[0,\frac{1}{d}\biggr), \\
\label{eqn:secondEq}
\frac{\partial v}{\partial t} &=  v (1-v) + \frac{\partial^{2} v}{\partial x^{2}} & \mbox{if} ~ v \in \biggl[\frac{1}{d},1\biggr].
\end{align}
The equation~\eqref{eqn:firstEq} can be rewritten to get
\begin{equation*}
\begin{split}
\frac1{\Delta x_i} \!\int_{Z_i} \frac{\partial v}{\partial t}(x,t)\,dx = & \,\frac{1}{\Delta x_i} \!\int_{Z_i} \!v(x,t)\bigl(1-v(x,t)\bigr)\,dx \\
&+ \frac{d}{\Delta x_i} \!\int_{Z_i} \frac{\partial }{\partial x}\left(v(x,t) \frac{\partial v}{\partial x}(x,t)\right) dx
\end{split}
\end{equation*}
where the finite volume integral average for the diffusion is dealt exactly as in the two-equations-based reduction, so that
\begin{equation*}
\begin{split}
& \frac{d}{\Delta x_i} \left[ v(x_{i+\frac12},t) \frac{\partial v}{\partial x}(x_{i+\frac12},t) - v(x_{i-\frac12},t) \frac{\partial v}{\partial x}(x_{i-\frac12},t) \right]\\
&\simeq \frac{d}{\Delta x_i} \left[ \frac{v_i(t)\Delta x_i+v_{i+1}(t)\Delta x_{i+1}}{\Delta x_i+\Delta x_{i+1}} \cdot \frac{v_{i+1}(t)-v_i(t)}{\dfrac{\Delta x_i}2+\dfrac{\Delta x_{i+1}}2} \right.\\
& \qquad\qquad \left. - \,\frac{v_{i-1}(t)\Delta x_{i-1}+v_i(t)\Delta x_i}{\Delta x_{i-1}+\Delta x_i} \cdot \frac{v_i(t)-v_{i-1}(t)}{\dfrac{\Delta x_{i-1}}2+\dfrac{\Delta x_i}2} \right].
\end{split}
\end{equation*}
Now, if the quantity $\Delta x_{i}$ is constant, we get the following semi-discrete version
\begin{equation*}
\begin{split}
\frac{d}{dt} v_i(t) = \,v_i(t)\bigl(1-v_i(t)\bigr) + \frac{d}{\Delta x} & \left[ \frac{v_i(t)\!+\!v_{i+1}(t)}{2} \cdot \frac{v_{i+1}(t)-v_i(t)}{\Delta x} \right.\\
& \;\left. - \,\frac{v_{i-1}(t)\!+\!v_i(t)}{2} \cdot \frac{v_i(t)-v_{i-1}(t)}{\Delta x} \right]
\end{split}
\end{equation*}
which easily leads to 
\begin{equation}
\label{eqn:firstsemiDiscr}
\begin{split}
\frac{d}{dt} v_i(t) &= \,v_i(t)\bigl(1-v_i(t)\bigr)  + \frac{d}{\Delta x^2}  \biggl[\frac{v_i(t)}{2} \bigl(v_{i+1}(t)-2\,v_i(t)+v_{i-1}(t)\bigr)\\
&\; + \frac{v_{i+1}(t)}{2} \bigl(v_{i+1}(t)-v_i(t)\bigr) - \frac{v_{i-1}(t)}{2} \bigl(v_i(t)-v_{i-1}(t)\bigr) \biggr].
\end{split}
\end{equation}
As concerns the equation~\eqref{eqn:secondEq}, by following the same path, in case of nonuniform mesh we have
\begin{equation*}
\frac{d}{dt} v_i(t) = \,v_i(t)\bigl(1-v_i(t)\bigr) + \frac1{\Delta x_i} \left[ \frac{v_{i+1}(t)-v_i(t)}{\dfrac{\Delta x_i}2+\dfrac{\Delta x_{i+1}}2} - \frac{v_i(t)-v_{i-1}(t)}{\dfrac{\Delta x_{i-1}}2+\dfrac{\Delta x_i}2} \right].
\end{equation*}
while, setting $\Delta x_{i}$ as a constant value, 
\begin{equation}
\label{eqn:secondsemiDiscr}
\frac{d}{dt} v_i(t) = \,v_i(t)\bigl(1-v_i(t)\bigr) + \frac{v_{i+1}(t)-2\,v_i(t)+v_{i-1}(t)}{\Delta x^2}\,.
\end{equation}
For the time discretization of ~\eqref{eqn:firstsemiDiscr} and \eqref{eqn:secondsemiDiscr}, we simply employ an explicit strategy, where $\Delta t\,$ is the fixed time step, so that the final numerical scheme reads as
\begin{equation}
\label{eqn:OneRedDiscr}
v^{n+1}_i =  v^n_ i +  \Delta t\,v^n_i \bigl(1-v^n_i \bigr) + H\bigl(v^{n}_ {i+1},v^{n}_ {i},v^{n}_ {i-1},\Delta x,\Delta t,d \bigr)
\end{equation}
where the function $H$ is defined as
\begin{equation*}
H= \begin{cases}
\begin{aligned}
\,d \dfrac{\Delta t}{\Delta x^{2}} \biggl[ \dfrac{v^{n}_{i}}{2} \bigl(v^{n}_{i+1}-2\,v^{n}_i+v^{n}_{i-1} \bigr)+\,\dfrac{v^{n}_{i+1}}{2} \bigl(v^{n}_{i+1}-v^{n}_i \bigr)-\,&\dfrac{v^{n}_{i-1}}{2} \bigl(v^{n}_i-v^{n}_{i-1} \bigr) \biggr],\\
&\mbox{if} ~  v^n_ i \in \biggl[0,\dfrac{1}{d}\biggr)  
\end{aligned}
& \\  \dfrac{\Delta t}{\Delta x^2} \bigl(v^{n}_ {i-1}-2\,v^{n}_i+v^{n}_{i+1} \bigr),\qquad  \mbox{if} ~  v^n_ i \in \biggl[\dfrac{1}{d},1\biggr].
\end{cases}
\end{equation*}
It is noticeable that the function $H$ exhibits a jump due to the discontinuity located in $v=\frac{1}{d}$ for the derivative of the diffusion $F$. 
\subsection{Simulations results}
We take advantage of the numerical scheme in~\eqref{eqn:OneRedDiscr} and perform numerical simulations in order to validate our one-equation-based reduction~\eqref{eqn:OneEqmodel} for the Gatenby-Gawlinski model. As regards the initial profile, we consider the Riemann problem whose states are $\mathbf{P_{L}}=1$ at the left and $\mathbf{P_{R}}=0$ at the right; all the parameters employed are listed in Table~\ref{tab:parametersOneEqReduc}.
\begin{table}[!ht]
\caption{Numerical default values for the parameters involved in the one-equation-based reduction.}
\label{tab:parametersOneEqReduc}
\begin{tabular}{*{4}{c}}
\toprule
$\mathbf{d}$ & $\mathbf{\Delta x}$ & $\mathbf{\Delta t}$  & $\mathbf{T}$\\
\midrule
$\{0.5,2\}$ & $0.05$ & $0.001$ & $20$\\
\bottomrule
\end{tabular}
\end{table} 

\begin{figure}[!ht]
\subfloat[][\emph{heterogeneous case}]
{\includegraphics[width=.49\columnwidth]{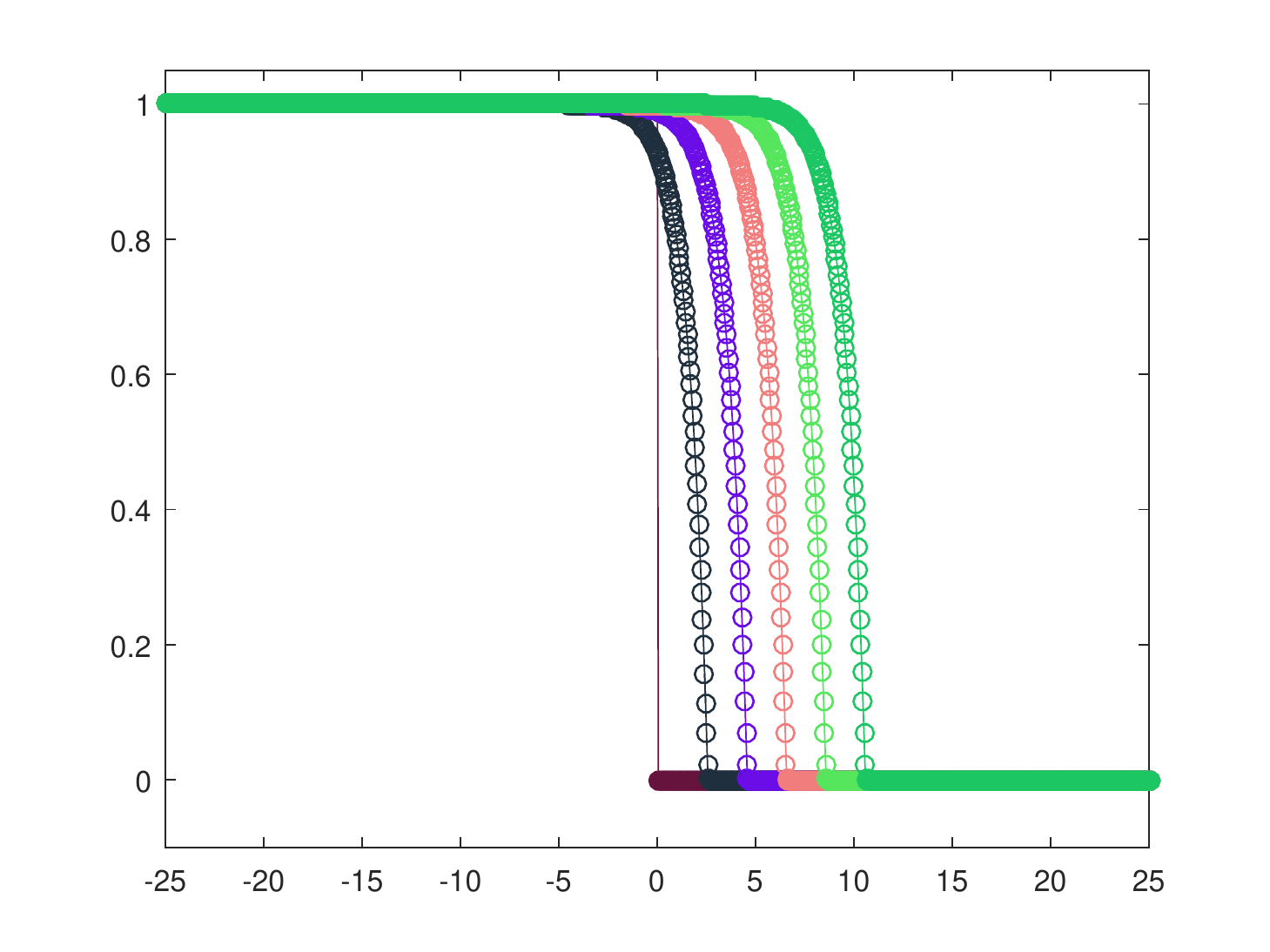}}
\subfloat[][\emph{homogeneous case}]
{\includegraphics[width=.49\columnwidth]{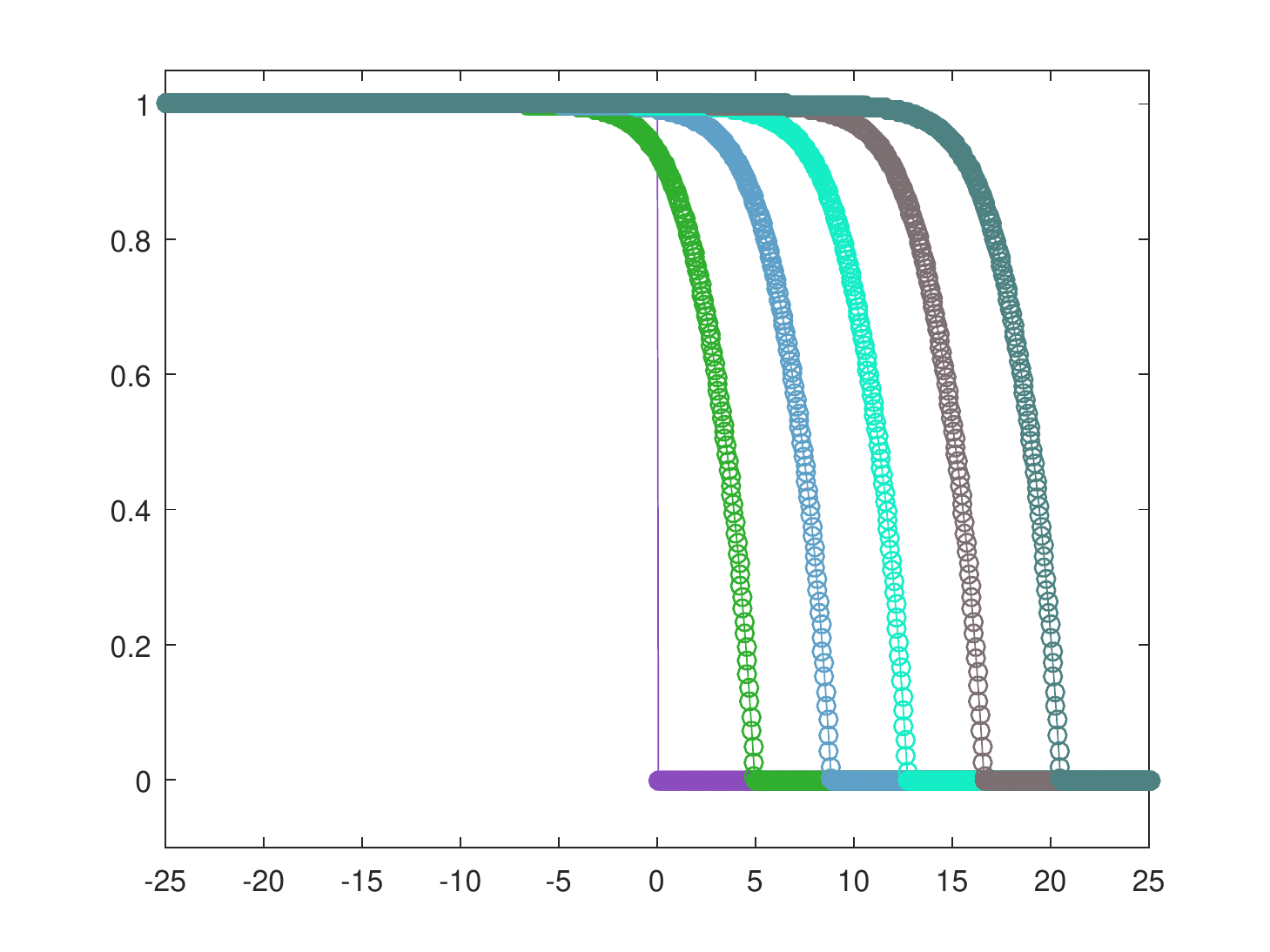}}
\caption{Front evolution for the tumour cells density in the heterogeneous case (A) and the homogeneous case (B). The parameters used are listed in Table~\ref{tab:parametersOneEqReduc}.}
\label{fig:tumourOneEq}
\end{figure}

The results are depicted in Figure~\ref{fig:tumourOneEq}(A) and Figure~\ref{fig:tumourOneEq}(B) for the heterogeneous and homogeneous invasion, respectively, by means of the front evolution representation. The cancerous cells density is plotted at equally spaced time instants and, in both the cases, the corresponding traveling waves turn out to be of sharp-type. 

\begin{figure}[!ht]
\includegraphics[width=.55\textwidth]{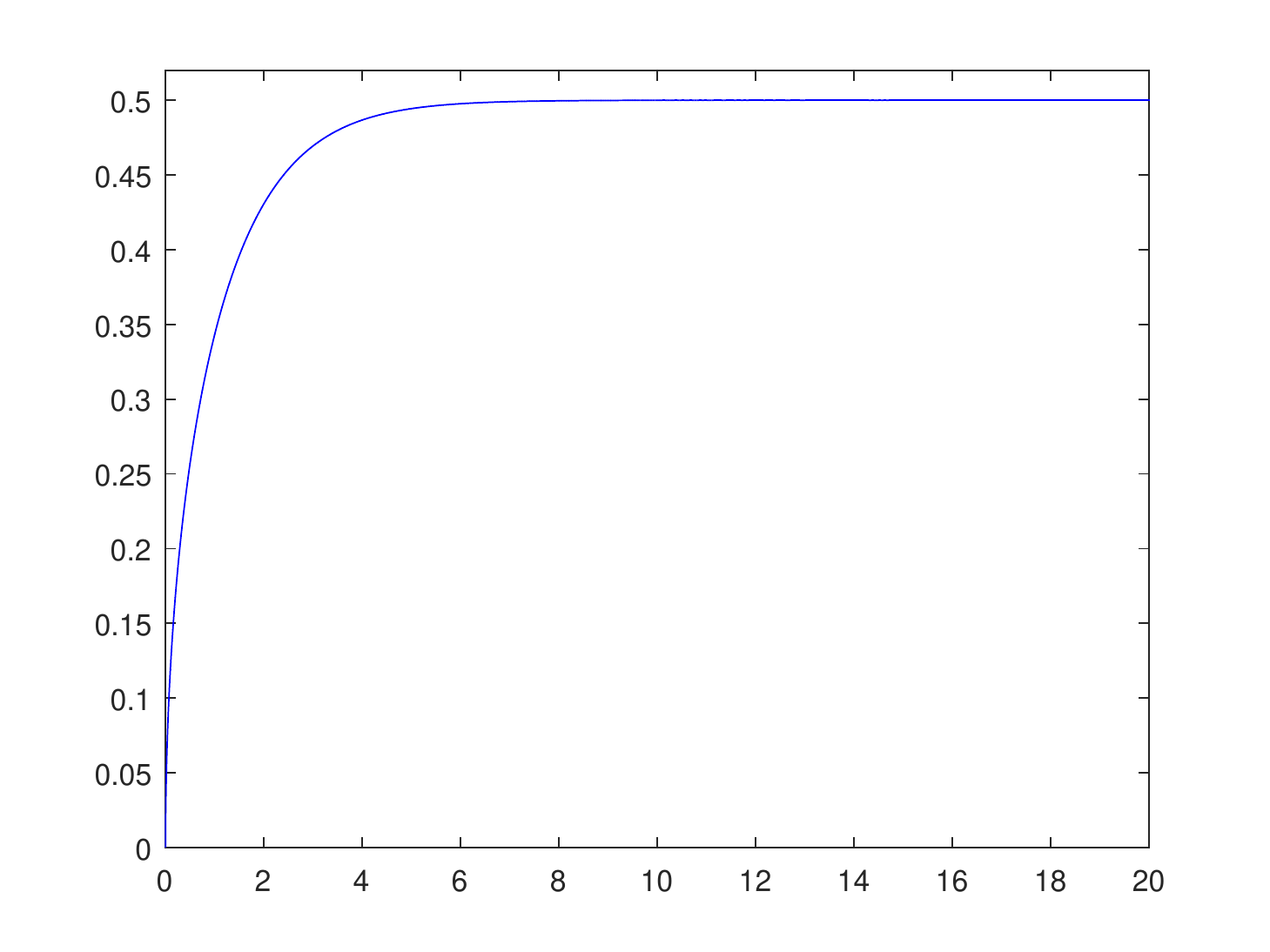}
\caption{Space-averaged propagation speed approximations as a function of time for the heterogeneous invasion. The parameters used are listed in Table~\ref{tab:parametersOneEqReduc}.}
\label{fig:speedheter}
\end{figure}

For ensuring the effectiveness of the one-equation-based reduction~\eqref{eqn:OneEqmodel}, it is important to establish if trends related to tumour invasions are correctly caught. In this respect, Figure~\ref{fig:tumourOneEq}, as well as providing information about the sharpness of the fronts, certifies as cancerous cells front moves forward faster in the homogeneous invasion. Specifically, adopting the space-averaged estimate~\eqref{eqn:numspeed}, we get $s \approx 0.499958$ for the heterogeneous case and $s \approx 0.968813$ for the homogeneous case: these two values are the asymptotic wave speeds of the tumour front. Figure~\ref{fig:speedheter} shows the discrete wave speed approximation computed as a function of time ($d=0.5$ is taken as a sample) and allows us to appreciate the convergence towards the corresponding asymptotic threshold.

\begin{figure}[!ht]
\subfloat[][\emph{heterogeneous invasion}]
{\includegraphics[width=.49\columnwidth]{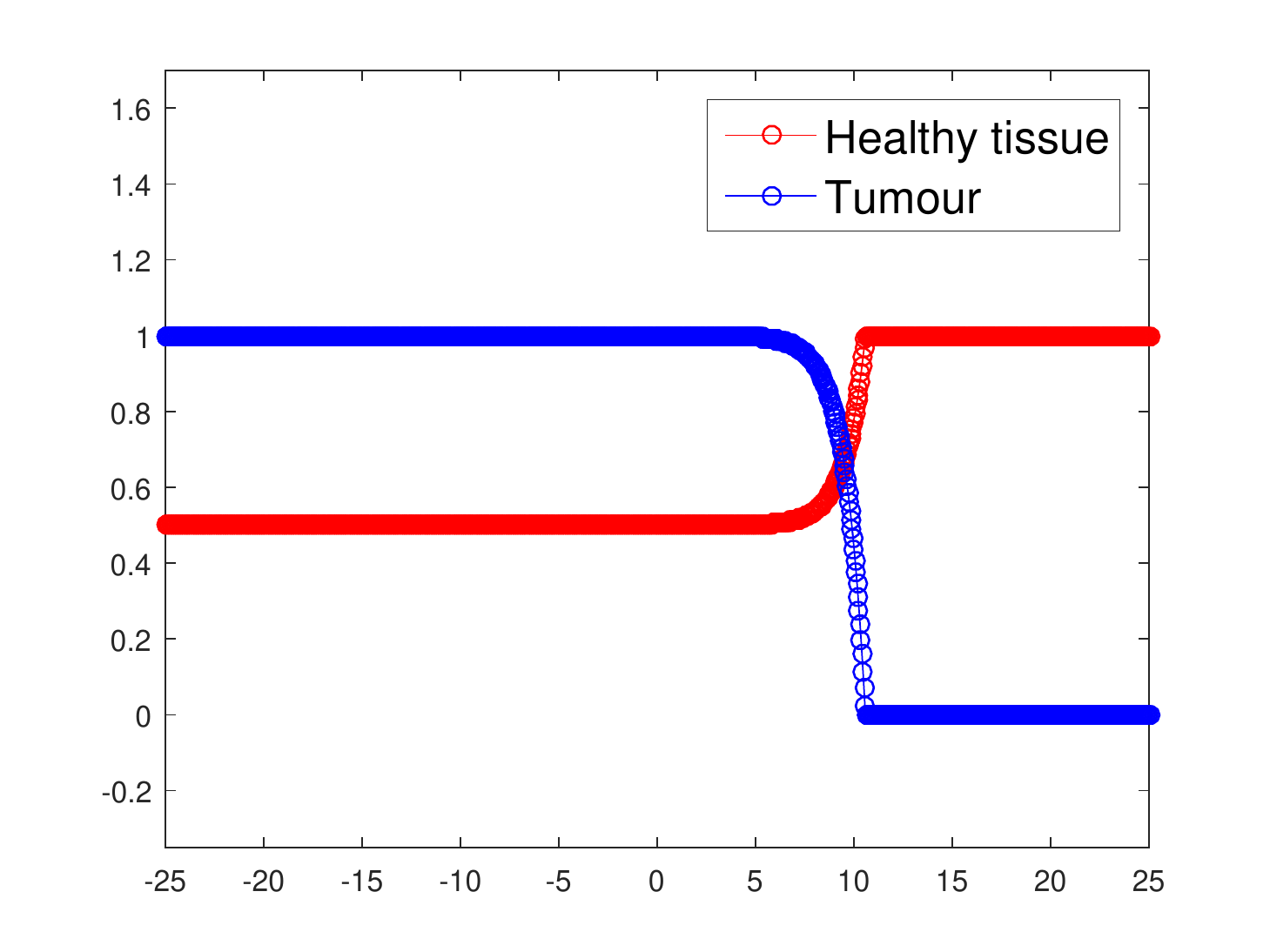}}
\subfloat[][\emph{homogeneous invasion}]
{\includegraphics[width=.49\columnwidth]{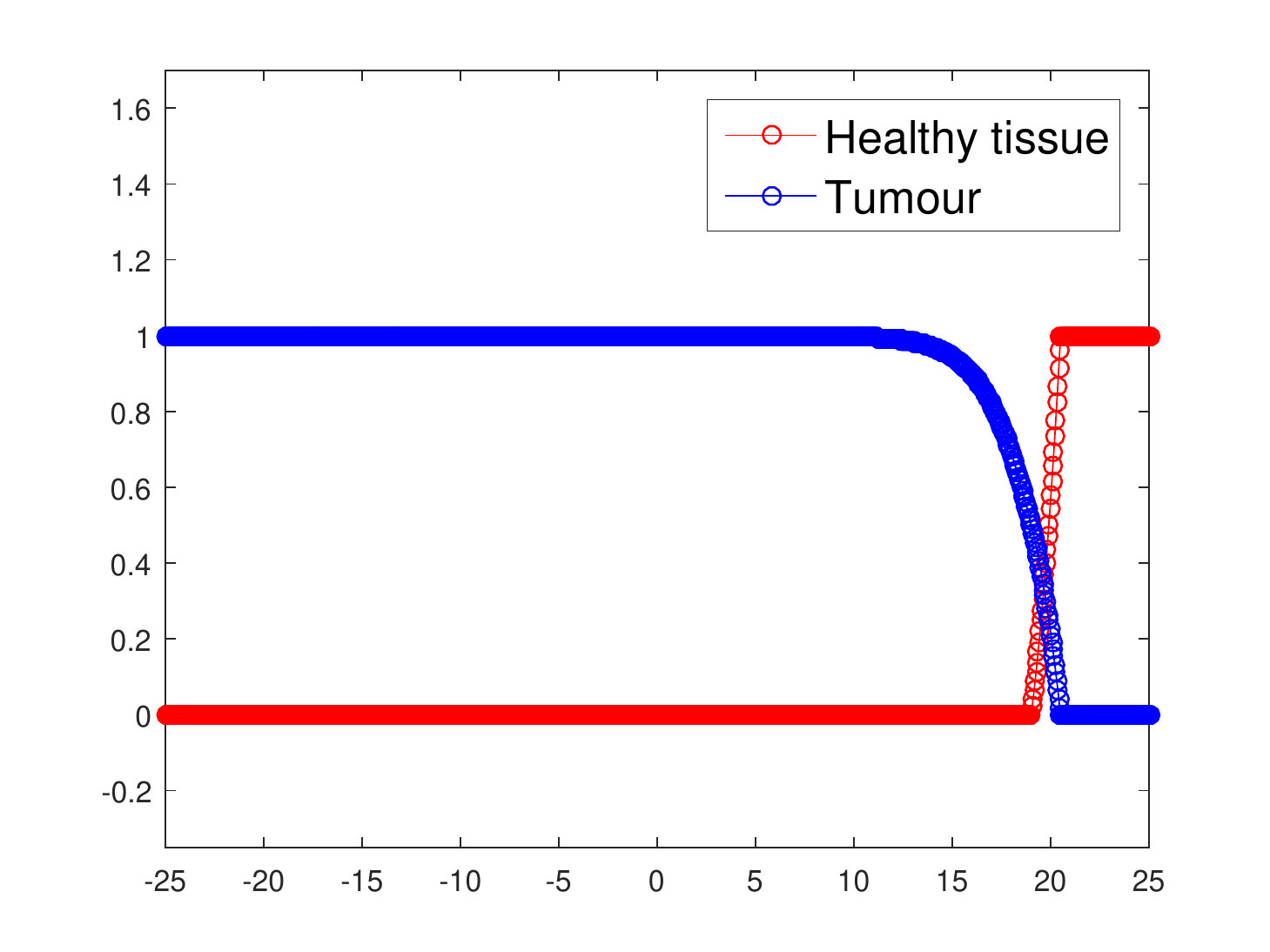}}
\caption{Numerical approximation of the tumour cells density along with the corresponding healthy cells profile, recovered by means of~\eqref{eqn:OneEqHp}, for both the heterogeneous case (A) and the homogeneous case (B). The parameters used are listed in Table~\ref{tab:parametersOneEqReduc}.}
\label{fig:solutionsOneEq}
\end{figure}

Finally, as further evidence of the reliability of~\eqref{eqn:OneEqmodel}, we can get information about the healthy cells density too, by simply leaning on~\eqref{eqn:OneEqHp}. The related graphs are depicted in Figure~\ref{fig:solutionsOneEq}(A) and Figure~\ref{fig:solutionsOneEq}(B). The plots are realized by simultaneously reporting numerical approximations for the cancerous cells densities, along with the induced healthy cells densities defined by means of~\eqref{eqn:OneEqHp}. The results qualitatively line up with the corresponding ones achieved for the two-equations-based reduction in Section~\ref{sec:twoeq}, so that the characteristic trends proper of the Gatenby-Gawlinski model are globally retrieved.

\begin{figure}[!ht]
\includegraphics[width=.55\textwidth]{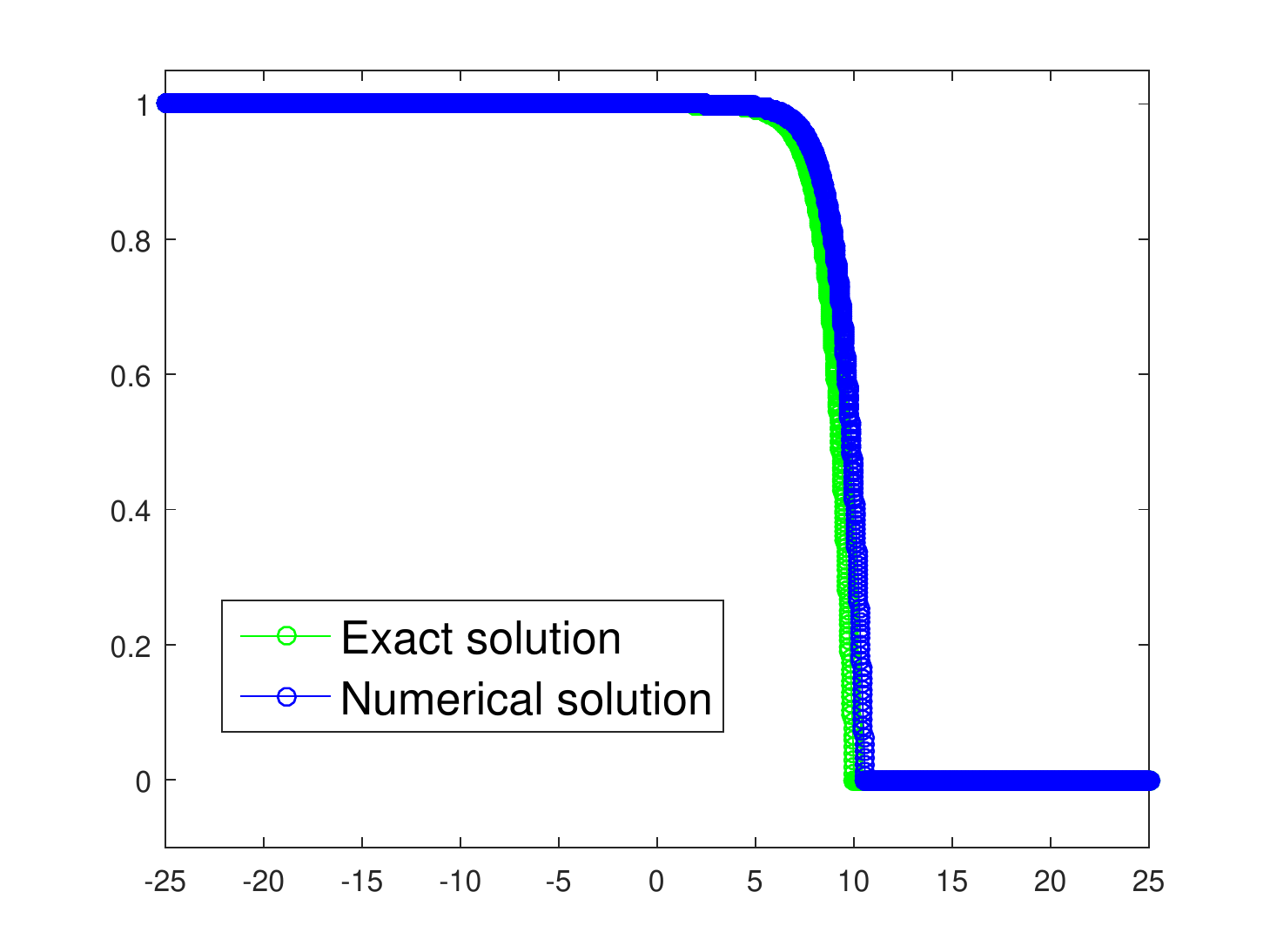}
\caption{Exact solution and corresponding tumour density numerical approximation for the one-equation-based model in case of heterogeneous invasion. $\Delta x = 0.01$ and $\Delta t = 0.0001$ are the choices for the spacial and time steps, respectively.}
\label{fig:OneEqExact}
\end{figure}

As far as the model described by~\eqref{eqn:OneEqmodel}, it is useful to notice that at least for the heterogeneous invasion, namely when $d<1$, the diffusion term turns out to be identically $F(v)=dv$: in order to easily check this statement, recalling the $F(v)$ definition in~\eqref{eqn:DegDiffPLF}, it follows that $v$ is always in $[0,1/d)$ if $d<1$, being $1/d>1$ and taking in mind the constraint $v \leq 1$. Due to this fact, the one-equation-based model in this specific case becomes
\begin{equation}
\label{eqn:OneEqmodelHeter}
\frac{\partial v}{\partial t} =  v (1-v) + d\frac{\partial}{\partial x}\left(v\frac{\partial v}{\partial x}\right).
\end{equation}   
For degenerate reaction-diffusion equations such as~\eqref{eqn:OneEqmodelHeter}, it is possible to get an analytical solution~\cite{Newman 1980, Newman 1983}. In this context, we simply impose that $v(x,t)$ is a propagating front of the form $\phi(x-st)$ being $s$ the associated wave speed and, after some conventional operations, the exact solution reads as
\begin{equation}
\label{eqn:exact}
v(x,t)= \begin{cases} 1-\exp\biggl(\dfrac{1}{\sqrt{2d}}(x-st)\biggr) &  \mbox{if} ~ x \leq st \\  0 & \mbox{if} ~ x > st,
\end{cases}
\end{equation}
where $s=\sqrt{d/2}$. Assuming the previous choice $d=0.5$, it follows that $s=0.5$, which is a threshold very close to our numerical estimate $s \approx 0.499958$. For the sake of completeness, we provide a graphical check as well, the plot being depicted in Figure~\ref{fig:OneEqExact}. We have chosen to exploit a refined spatio-temporal mesh, namely $\Delta x = 0.01$ and $\Delta t = 0.0001$, compared to the parameters listed in Table~\ref{tab:parametersOneEqReduc}, with the aim of achieving a very effective graphical result. The resulting trajectories are very close and the wave speed approximation is good too, being $s \approx 0.499983$.
  
Finally, by analogy with what has been shown in~\cite{Moschetta 2019} in order to numerically appreciate the transition occurring from the complete Gatenby-Gawlinski model towards the two-equations-based reduction by increasing the parameter $c$ in~\eqref{eqn:system}, we propose a similar analysis regarding the one-equation-based model. Indeed, recalling the assumption~\eqref{eqn:OneEqHp} exploited in~\eqref{eqn:system3} for justifying the model simplification, it is possible studying the transition occurring between the two-equations-based and one-equation-based model by defining the $\epsilon$-dependent time derivative of the function $u(x,t)$. We get
\begin{equation}
\label{eqn:system4}
\begin{dcases}
\epsilon \frac{\partial u}{\partial t} = u (1-u) - d u v\\
\frac{\partial v}{\partial t} = v (1-v) + \frac{\partial}{\partial x} \left[(1-u) \frac{\partial v}{\partial x}\right].
\end{dcases}
\end{equation}
At this stage, we can easily infer that, taking the limit as the parameter $\epsilon$ approaches zero in the first equation of~\eqref{eqn:system4}, perfectly matches, from a theoretical point of view, the idea behind the hypotesis~\eqref{eqn:OneEqHp}, which automatically leads to~\eqref{eqn:OneEqmodel}. Now, we want to catch the transition, either employing the wave speed numerical estimate~\eqref{eqn:numspeed}, which is the approach proposed in~\cite{Moschetta 2019}, either taking advantage of the solution~\eqref{eqn:exact}. As a matter of fact, considering the heterogeneous invasion context, we can rely on the exact solution for the one-equation-based model and exploit it to verify the transition from the two-equations reduction. 
\begin{figure}[!ht]
\subfloat[][\emph{convergence of the wave speeds estimates}]
{\includegraphics[width=.49\columnwidth]{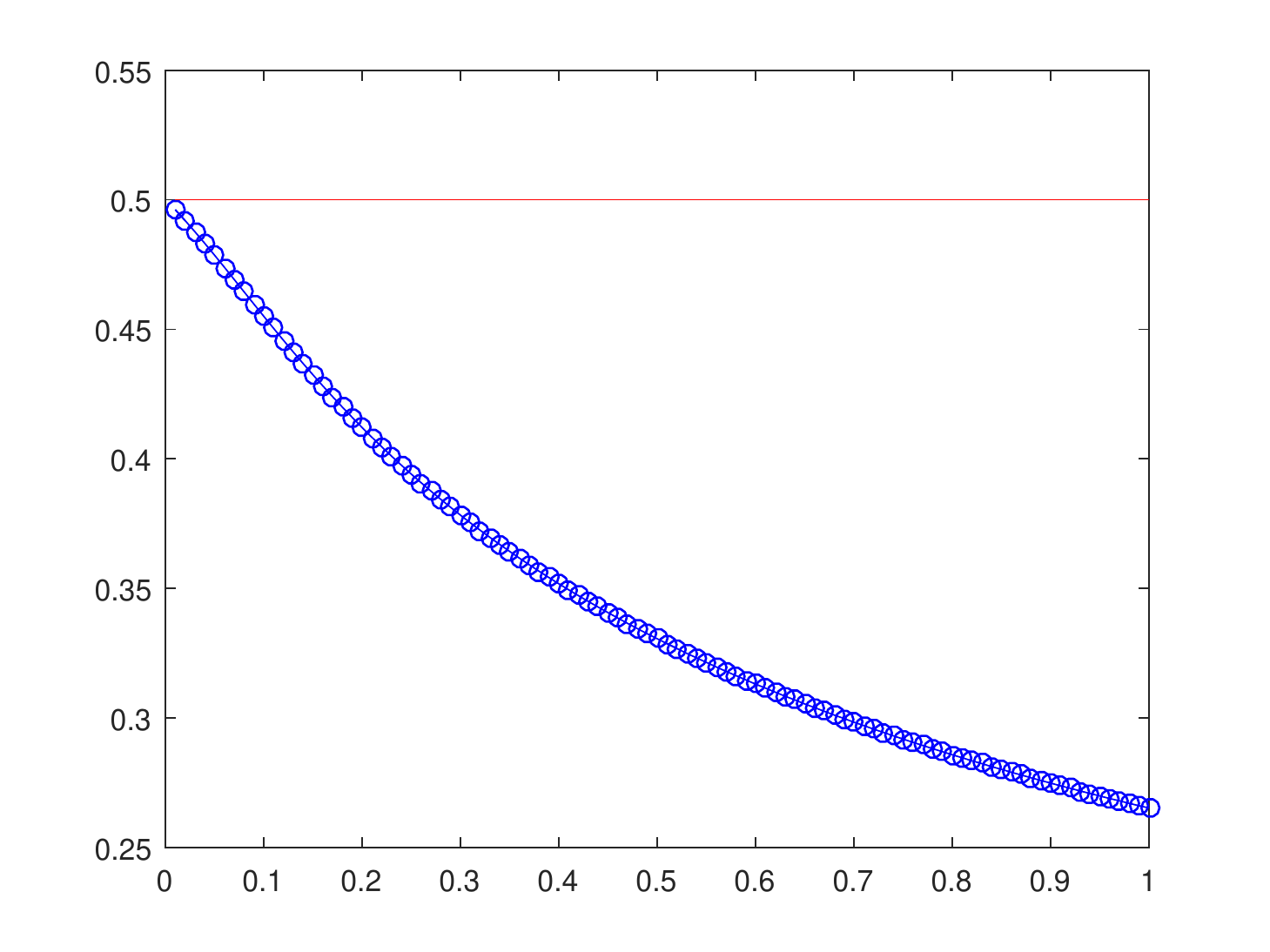}}
\subfloat[][\emph{convergence of the tumour density numerical solution}]
{\includegraphics[width=.49\columnwidth]{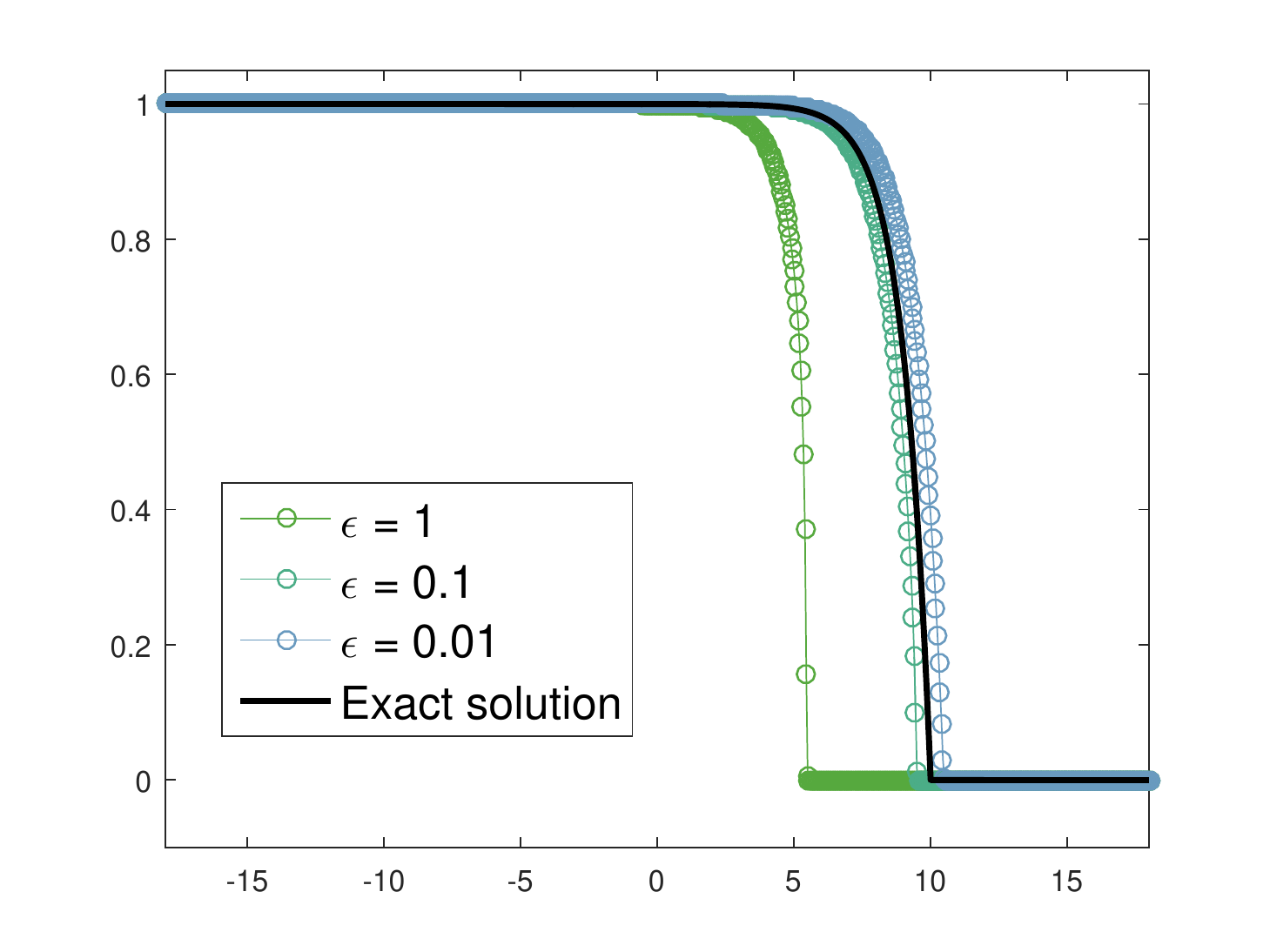}}
\caption{Wave speeds estimates (blue circles) for $v(x,t)$ in~\eqref{eqn:system4} if $\epsilon \in [0.01,1]$ along with the wave speed analytical value (red line) provided by~\eqref{eqn:exact} (A) and convergence of the tumour density numerical approximation from~\eqref{eqn:system4}, as a function of specific $\epsilon$ values, towards the exact solution~\eqref{eqn:exact} (B). The parameters used are listed in Table~\ref{tab:parametersOneEqReduc}.}
\label{fig:Eps}
\end{figure}
Figure~\ref{fig:Eps}(A) exhibits the wave speed numerical approximations achieved by assuming decreasing $\epsilon$ values in~\eqref{eqn:system4} in the case of tumour cells front. The resulting trend correctly reports a convergence towards the asymptotic threshold $s=0.5$, which is the analytical prediction for~\eqref{eqn:OneEqmodelHeter}. Moreover, the exact solution~\eqref{eqn:exact} allows us to graphically recognize the transition towards~\eqref{eqn:OneEqmodelHeter} by means of a convergence check: Figure~\ref{fig:Eps}(B) reports the tumour density numerical approximation provided by~\eqref{eqn:system4} as a function of some $\epsilon$ values taken as sample. It is possible to detect a progressive alignment with the analytical solution~\eqref{eqn:exact}.

%%%------------------------------------------------------------------------------------------
\section{Conclusions}
\label{sec:theend}
In this paper, we have investigated the Gatenby-Gawlinski model for tumour invasion taking advantage of some useful system simplifications. In continuity with the study carried out in~\cite{Moschetta 2019}, we have extended the available results for the two-equations-based reduction, specifically analyzing the sharpness/smoothness of the fronts by means of numerical simulations based on a finite volume approximation. A qualitative check about the shape of the traveling waves suggest that we deal with smooth-type fronts concerning the full model and with sharp-type ones for the two-equations reduction. A sensitivity analysis with respect to the system parameters $r$ and $d$ in~\eqref{eqn:system3} has been provided as well, by using as unknown the wave speed, whose numerical approximation is achieved through a space-averaged estimate~\cite{LeVeque 1990, Moschetta 2019}.

Subsequently, we have proposed a further system simplification leading to a one-equation-based reduction, framed within the degenerate reaction-diffusion equations field. Several results are available in the literature about the existence and uniqueness of sharp/smooth-type fronts for such a mathematical problem~\cite{Malaguti 2003, Sanchez 1994 first, Sanchez 1995, Sanchez 1996}, but as regards our case~\eqref{eqn:OneEqmodel}, due to the almost everywhere differentiability in $[0,1]$ of the degenerate diffusion, relying on numerical checks is required in order to establish the fronts shape.
The evidence is that traveling wave arising from our one-equation reduction proves itself to be of sharp-type. We have shown the reduction to qualitatively catch the typical dynamics of the Gatenby-Gawlinski model and, in the specific case of heterogeneous invasion, we have provided the corresponding analytical solution as well, exploiting its availability to verify the simulations results effectiveness. Finally, we have studied the transition from the two-equations-based reduction towards the one-equation-based simplification by defining the $\epsilon$-dependent time derivative of the healthy cells density $u(x,t)$, as in~\eqref{eqn:system4}.

As far as the necessity to lean on numerical assessments about the sharpness detection in the one-equation reduction case, arising from the requirement for a more regular degenerate diffusion in terms of the available theoretical results~\cite{Malaguti 2003, Sanchez 1995}, the possibility of employing a smoother approximation to by-pass the almost everywhere differentiability of the diffusion, seems to be a very promising ground. In this way, enough regularity would be ensured to prove existence and uniqueness results for the traveling fronts.

%%%------------------------------------------------------------------------------------------

\end{document}